\documentclass[final,1p,times]{elsarticle}%
\usepackage{lineno,hyperref}
\usepackage{graphicx}
\usepackage{amsmath}
\usepackage{amsfonts}
\usepackage{amssymb}
\usepackage{graphicx}
\usepackage{multirow}
\usepackage{enumitem}
\usepackage{color}%

\usepackage{amsmath,multirow}
\usepackage{amssymb}
\usepackage{amscd}
\usepackage{amsfonts}
\usepackage{enumerate}
\usepackage{mathrsfs}
\usepackage{xypic}
\usepackage{stmaryrd}
\usepackage{graphicx}
\usepackage{fancybox}
\usepackage{portland}
\usepackage{color}
\usepackage{amsmath}
\usepackage{exscale}
\usepackage{enumitem,array}%
\usepackage{amssymb}
\newtheorem{mytheo}{Theorem}[section]
\newtheorem{mydef}[mytheo]{Definition}

\newtheorem{lem}[mytheo]{Lemma}
\newtheorem{rem}{Remark}
\newcommand{\norm}[1]{\left\Vert#1\right\Vert}
\journal{Journal of  Computational and Applied Mathematics}

\begin{document}

\begin{frontmatter}

\title{Trigonometric collocation methods based on Lagrange basis polynomials for
multi-frequency oscillatory second-order differential
equations\tnoteref{mytitlenote}} \tnotetext[mytitlenote]{This paper
was supported by  National Natural Science Foundation of China under
Grant 11401333,11271186,11171178, by  Natural Science Foundation of
Shandong Province under Grant ZR2014AQ003, and by  China
Postdoctoral Science Foundation under Grant 2015M580578. }

\author[Wang]{Bin~Wang\corref{cor1}}
\author[Wu]{Xinyuan Wu}
\author[meng]{Fanwei Meng}
\address[Wang]{School of Mathematical Sciences, Qufu Normal University, Qufu, Shandong 273165, P.R.China}
\address[Wu]{Department of Mathematics, Nanjing University, Nanjing 210093, P.R.China}
\address[meng]{School of Mathematical Sciences, Qufu Normal University, Qufu, Shandong 273165, P.R.China}

 \ead{wangbinmaths@gmail.com (Bin Wang), xywu@nju.edu.cn (Xinyuan Wu), fwmeng@mail.qfnu.edu.cn}
\cortext[cor1]{Corresponding author.}

\begin{abstract}
In the present work, a kind of trigonometric collocation methods
based on Lagrange basis polynomials  is developed  for effectively
solving multi-frequency oscillatory second-order differential
equations $q^{\prime\prime}(t)+Mq(t)=f\big(q(t)\big)$. The
properties of the obtained methods are investigated. It is shown
that the convergent condition   of these   methods is independent of
$\norm{M}$, which is very  crucial  for solving oscillatory systems.
A fourth-order scheme  of the  methods is presented. Numerical
experiments are implemented to show the remarkable efficiency of the
methods proposed in this paper.
\end{abstract}


\end{frontmatter}

\section{Introduction}
The  numerical treatment  of multi-frequency oscillatory systems is
a computational problem of an overarching importance in a wide range
of applications,  such as  quantum physics, circuit simulations,
flexible body dynamics and mechanics  (see, e.g.
\cite{cohen2005bit,Lubich2006,B.garcia1999,hairer2000,hairer2006,wu2012-6,wu2013-book}
 and the references therein). The  main theme   of the present paper is to  construct and analyse a kind of   efficient collocation methods for
solving multi-frequency oscillatory second-order differential
equations of the form
\begin{equation}
q^{\prime\prime}(t)+Mq(t)=f\big(q(t)\big),  \qquad
q(0)=q_0,\ \ q'(0)=q'_0,\qquad t\in[0,t_{\mathrm{end}}],\label{prob}%
\end{equation}
where $M$ is a $d\times d$ positive semi-definite matrix
 implicitly containing the  frequencies of the
oscillatory problem and $f: \mathbb{R}^{d}\rightarrow
\mathbb{R}^{d}$ is  an analytic function. The solution of this
system is a multi-frequency nonlinear oscillator because of the
presence of the linear term $Mq$. System \eqref{prob} is  a highly
oscillatory problem when  $\norm{M}\gg1$. In recent years, various
numerical methods for   approximating solutions of oscillatory
systems
  have been developed by many researchers.  Readers are referred
  to \cite{cohen2005bit,Hochbruck1999,ANM-2015,JCP-2013,wang2012-1,wu2012-5,wu2013-ANM,wu2012-6,wu-2012-BIT,wu2013-book}
and the references therein. Once it  is further assumed that  $M$ is
symmetric and $f$ is the negative gradient of a real-valued function
$U(q)$, the system \eqref{prob} is identical to the following
initial value Hamiltonian system
\begin{equation}
\begin{aligned}
 &\dot{q}=\nabla_p H(q,p),\qquad \ \ q(0)=q_{0},\\
 &\dot{p}=-\nabla_q H(q,p),\qquad  p(0)=p_{0}\equiv q'_0 \end{aligned}
\label{H-s}%
\end{equation}
with the Hamiltonian function
\begin{equation}H(q,p)=\frac{1}{2}p^{\intercal}p+\frac{1}{2}q^{\intercal}Mq+U(q).\label{H}%
\end{equation}
This is an  important  system which  has  received much attention by
many authors (see, e.g.
\cite{Cohen-2006,cohen2005bit,Lubich2006,hairer2000,hairer2006}).

In \cite{wang-2014},  the authors  took advantage of shifted
Legendre polynomials to obtain a local Fourier expansion of the
system \eqref{prob} and derived a  kind of collocation methods
(trigonometric collocation methods).
   The analysis and the results of numerical experiments in \cite{wang-2014} showed  that the
trigonometric  collocation methods are
 more efficient in comparison with some alternative approaches that
have previously appeared in the literature. Motivated by the work in
\cite{wang-2014},  this paper is devoted to the formulation and
analysis of another trigonometric collocation methods for solving
multi-frequency oscillatory second-order systems \eqref{prob}. We
will consider a more classical approach and use Lagrange polynomials
to obtain the methods. Because of this different approach, compared
with the methods in \cite{wang-2014}, the obtained methods  have a
simpler scheme and can  be implemented in practical computations at
a lower cost. These trigonometric collocation
  methods  are  designed  by interpolating the function $f$ of \eqref{prob} by Lagrange basis
polynomials,  and incorporating the variation-of-constants formula
with  the idea of collocation methods. It is noted that  these
integrators  are a kind of collocation methods and they share all
the interesting features of collocation methods. We  analyse the
properties of the
 trigonometric collocation methods. We also consider  the convergence of the fixed-point
iteration for the    methods.    It is  important to emphasize that
for the   trigonometric collocation methods, the convergent
condition is independent of $\norm{M}$, which is a very
 important property for solving oscillatory systems.

 This paper is organized as follows. In Section \ref{sec:Formulation}, we formulate the scheme
  of trigonometric collocation methods based on Lagrange
basis polynomials.   The properties of the obtained methods are
analysed in Section \ref{sec:Analysis of the methods}. In Section
\ref{numerical experiments}, a fourth-order scheme  of the  methods
is presented and numerical tests confirm that the  method proposed
in this paper yields a dramatic improvement. Conclusions are
included in Section \ref{sec:conclusion}.

\section{Formulation of the   methods} \label{sec:Formulation}
To begin with we  restrict the multi-frequency oscillatory system
\eqref{prob}  to the interval $[0,h]$ with any $h>0$:
\begin{equation}
q''(t)+Mq(t)=f\big(q(t)\big),  \qquad
q(0)=q_0,\ \ q'(0)=q'_0,\qquad t\in[0,h].\label{prob h}%
\end{equation}
  With regard to  the variation-of-constants formula for
\eqref{prob} given in \cite{wu2010-1}, we have the following result
on the exact solution $q(t)$ of the system~\eqref{prob} and its
derivative $q'(t)=p(t)$:
  \begin{equation}
\begin{aligned} &q(t)=\phi_0(t^2M)q_0+t\phi_1(t^2M)p_0+ t^2\int_{0}^1(1-z)\phi_1\big((1 -z)^2t^2M\big)f\big(q(tz)\big)dz,\\
&p(t)=-tM\phi_1( t^2M)q_0+\phi_0(t^2M)p_0
+t\int_{0}^{1}\phi_0\big((1 -z)^2t^2M\big)f\big(q(tz)\big)dz,
\end{aligned}
\label{probsystemF3}%
\end{equation}
 where  $t\in[0,h]$  and \begin{equation}
\phi_{i}(M):=\sum\limits_{l=0}^{\infty}\frac{(-1)^{l}M^{l}}{(2l+i)!},\qquad
\ i=0,1. \label{Phi01}
\end{equation}
From this result, it follows that
  \begin{equation}
\begin{aligned} &q(h)=\phi_0(V)q_0+h\phi_1(V)p_0+ h^2\int_{0}^1(1-z)\phi_1\big((1 -z)^2V\big)f\big(q(hz)\big)dz,\\
&p(h)=-hM\phi_1( V)q_0+\phi_0(V)p_0 +h\int_{0}^{1}\phi_0\big((1
-z)^2V\big)f\big(q(hz)\big)dz,
\end{aligned}
\label{qh-ph}%
\end{equation}
 where   $V=h^2M.$

The main point in designing practical schemes to solve \eqref{prob}
is based on  replacing $f(q)$ in \eqref{qh-ph}
  by some expansion. In this paper, we interpolate $f(q)$  as
\begin{equation}f\big(q(\xi h)\big)\sim\sum\limits_{j=1}^
{s}l_j(\xi )f\big(q(c_j h)\big),\qquad   \xi \in[0,1],\label{fq}%
\end{equation}
where
\begin{equation}l_j(x)=\prod\limits_{k=1,k\neq j}^
{s}\frac{x-c_k}{c_j-c_k}\label{lj}%
\end{equation}
for $j=1,\ldots,s$ are the Lagrange basis polynomials in
interpolation and $c_1, \ldots , c_s$ are distinct real numbers
(usually $s\geq1,\ 0 \leq c_i \leq 1$). Then replacing $f(q(\xi h))$
in \eqref{qh-ph}
  by the series \eqref{fq} yields an approximation of $q(h),p(h)$    as follows:
\begin{equation}
\begin{aligned} &\tilde{q}(h)=\phi_0(V)q_0+h\phi_1(V)p_0+
h^2\sum\limits_{j=1}^
{s}I_{1,j}f\big(\tilde{q}(c_j h)\big),\\
&\tilde{p}(h)=-hM\phi_1( V)q_0+\phi_0(V)p_0 +h\sum\limits_{j=1}^
{s}I_{2,j}f\big(\tilde{q}(c_j h)\big),
\end{aligned}
\label{basic method}%
\end{equation}
where \begin{equation}
\begin{aligned} &I_{1,j}:=\int_{0}^1l_j(z)(1-z)\phi_1\big((1 -z)^2V\big)dz,\ \
I_{2,j}:=\int_{0}^{1}l_j(z)\phi_0\big((1 -z)^2V\big)dz.
\end{aligned}
\label{I12j}%
\end{equation}
According to the  variation-of-constants formula
\eqref{probsystemF3} for \eqref{prob h}, the approximation
\eqref{basic method} satisfies the following system
\begin{equation}
\begin{aligned}
 & \tilde{q}'(\xi h)=\tilde{p}(\xi h), \qquad \qquad \qquad \qquad \qquad \ \   \tilde{q}(0)=q_{0},\\
 &\tilde{p}'(\xi h)=-M\tilde{q}(\xi h)+\sum\limits_{j=1}^
{s}l_j(\xi )f\big(\tilde{q}(c_j h)\big),\ \ \ \
\tilde{p}(0)=p_{0}.\end{aligned}
\label{new system}%
\end{equation}

In what follows we first approximate  $f\big(\tilde{q}(c_j h)\big),\
I_{1,j},\ I_{2,j}$ appearing in \eqref{basic method} and then   a
  kind of collocation methods can be formulated.

\subsection{The computation of $f(\tilde{q}(c_j h))$} \label{subsec:Discretization}

 It follows  from \eqref{new system}  that $\tilde{q}(c_i h),\
i=1,2,\ldots,s,$ can be obtained by solving the following discrete
problems:
\begin{equation}
\begin{aligned}
&\tilde{q}''(c_i h)+M\tilde{q}(c_i h)=\sum\limits_{j=1}^
{s}l_j(c_i )f\big(\tilde{q}(c_j h)\big),\ \ \  \tilde{q}(0)=q_{0},\ \tilde{q}'(0)=p_{0}.\\
\end{aligned}
\label{discrete quadrature truncating H-s}%
\end{equation}
By setting $\tilde{q}_i=\tilde{q}(c_i h)$ with $i=1,2,\ldots,s,$
\eqref{discrete quadrature truncating H-s} can be solved by the
variation-of-constants formula \eqref{probsystemF3} in the form:
\begin{equation*}
\begin{aligned}
\tilde{q}_i=&\phi_0(c_i^2V)q_0+c_ih\phi_1(c_i^2V)p_0+
(c_ih)^2\sum\limits_{j=1}^ {s}\tilde{I}_{c_i,j}f(\tilde{q}_j ),\quad
i=1,2,\ldots,s,
\end{aligned}
\label{discrete quadrature method}%
\end{equation*}
where \begin{equation}
\begin{aligned} &\tilde{I}_{c_i,j}:=\int_{0}^1l_j(c_iz)(1-z)\phi_1\big((1
-z)^2c_i^2V\big)dz,\qquad i, j=1,\ldots,s.
\end{aligned}
\label{Ijci}%
\end{equation}
\subsection{The computation of $I_{1,j},\ I_{2,j},\ \tilde{I}_{c_i,j}$} \label{subsec:Computations of I}
With the definition \eqref{lj}, the integrals $I_{1,j},\ I_{2,j},\
\tilde{I}_{c_i,j}$ appearing above can be computed by
\begin{equation*}
\begin{aligned} I_{1,j}=&\int_{0}^1l_j(z)(1-z)\phi_1\big((1
-z)^2V\big)dz\\
=&\prod\limits_{k=1,k\neq j}^ {s}\sum\limits_{l=0}^{\infty}\int_{0}^1\frac{z-c_k}{c_j-c_k}(1-z)^{2l+1}dz\frac{(-1)^{l}V^{l}}{(2l+1)!}\\
=&\sum\limits_{l=0}^{\infty}\Big(\prod\limits_{k=1,k\neq j}^
{s}\frac{\frac{1}{2l+3}-c_k}{c_j-c_k}\Big)\frac{(-1)^{l}V^{l}}{(2l+2)!}
=\sum\limits_{l=0}^{\infty}l_j\Big(\frac{1}{2l+3}\Big)\frac{(-1)^{l}V^{l}}{(2l+2)!},\\
 I_{2,j}=&\int_{0}^1l_j(z)\phi_0\big((1
-z)^2V\big)dz
=\prod\limits_{k=1,k\neq j}^ {s}\sum\limits_{l=0}^{\infty}\int_{0}^1\frac{z-c_k}{c_j-c_k}(1-z)^{2l}dz\frac{(-1)^{l}V^{l}}{(2l)!}\\
=&\sum\limits_{l=0}^{\infty}\Big(\prod\limits_{k=1,k\neq j}^
{s}\frac{\frac{1}{2l+2}-c_k}{c_j-c_k}\Big)\frac{(-1)^{l}V^{l}}{(2l+1)!}
=\sum\limits_{l=0}^{\infty}l_j\Big(\frac{1}{2l+2}\Big)\frac{(-1)^{l}V^{l}}{(2l+1)!},\\
\tilde{I}_{c_i,j}=&\int_{0}^1l_j(c_iz)(1-z)\phi_1\big((1
-z)^2c_i^2V\big)dz\\
=&\prod\limits_{k=1,k\neq j}^ {s}\sum\limits_{l=0}^{\infty}\int_{0}^1\frac{c_iz-c_k}{c_j-c_k}(1-z)^{2l+1}dz\frac{(-1)^{l}(c_i^2V)^{l}}{(2l+1)!}\\
=&\sum\limits_{l=0}^{\infty}\Big(\prod\limits_{k=1,k\neq j}^
{s}\frac{\frac{c_i}{2l+3}-c_k}{c_j-c_k}\Big)\frac{(-1)^{l}(c_i^2V)^{l}}{(2l+2)!}
=\sum\limits_{l=0}^{\infty}l_j\Big(\frac{c_i}{2l+3}\Big)\frac{(-1)^{l}(c_i^2V)^{l}}{(2l+2)!},\\
&\qquad \qquad\qquad\qquad\qquad\qquad\qquad \qquad\qquad\qquad
\qquad i,j=1,\ldots,s.
\end{aligned}
\label{Computate I}%
\end{equation*}
When the matrix $M$ is  symmetric and positive semi-definite, we
have the decomposition of $M$  as follows:
\[
M=P^{\intercal}W^{2}P=\Omega_{0}^{2}\ \ \mbox{with}\
\Omega_{0}=P^{\intercal}W P,
\]
where ${P}$ is an orthogonal matrix and $W=\textmd{diag}(\lambda_k)$
  with nonnegative diagonal entries which are the
square roots of the eigenvalues of $M$. Then the above  integrals
become

\begin{equation*}
\begin{aligned} I_{1,j}=&P^{\intercal}\int_{0}^1 l_j(z)  W ^{-1}\sin\big((1
-z)W\big)dzP=P^{\intercal}\textmd{diag}\Big(\int_{0}^1 l_j(z)
\lambda_k ^{-1}\sin\big((1
-z)\lambda_k\big)dz\Big)P,\\
 I_{2,j}=&P^{\intercal}\int_{0}^1l_j(z)\cos\big((1 -z)W\big)dzP=P^{\intercal}\textmd{diag}\Big(\int_{0}^1l_j(z)\cos\big((1
 -z)\lambda_k\big)dz\Big)P,\\
\tilde{I}_{c_i,j}=&P^{\intercal}\int_{0}^1 l_j(c_iz) (c_iW
)^{-1}\sin\big((1
-z)c_iW\big)dz=P^{\intercal}\textmd{diag}\Big(\int_{0}^1 l_j(c_iz)
(c_i\lambda_k )^{-1}\sin\big((1 -z)c_i\lambda_k\big)dz\Big)P,
\\
  & i,j=1,\ldots,s.
\end{aligned}
\end{equation*}
 It is noted that  $W^{-1}\sin\big((1
-z)W\big),\ (c_iW )^{-1}\sin\big((1 -z)c_iW\big)$ are well-defined
also for singular $W$.
 The case of $\lambda_k=0$ gives:
\begin{equation*}
\begin{aligned} \int_{0}^1 l_j(z)  \lambda_k
^{-1}\sin\big((1 -z)\lambda_k\big)dz=&\int_{0}^1 l_j(z)  (1
-z)dz,\\
\int_{0}^1l_j(z)\cos\big((1 -z)\lambda_k\big)dz=&\int_{0}^1l_j(z)dz,\\
\int_{0}^1 l_j(c_iz) (c_i\lambda_k )^{-1}\sin\big((1
-z)c_i\lambda_k\big)dz=&\int_{0}^1 l_j(c_iz) (1 -z)dz,
\end{aligned}
\end{equation*}
 which  can be evaluated easily since $l_j(z)$ is a
 polynomial function. If $\lambda_k\neq0$, they can be
evaluated as follows:\\ $
\begin{aligned}&\int_{0}^1 l_j(z)  \lambda_k
^{-1}\sin\big((1 -z)\lambda_k\big)dz \\
=&1/\lambda_k\int_{0}^1 l_{j}(z)\sin\big((1
-z)\lambda_k\big)dz=1/\lambda_k^2\int_{0}^1
l_{j}(z)d\cos\big((1 -z)\lambda_k\big)\\
=& 1/\lambda_k^2 l_{j}(1)-1/\lambda_k^2 l_{j}(0)\cos(\lambda_k)-
1/\lambda_k^2\int_{0}^1
l'_{j}(z)\cos\big((1 -z)\lambda_k\big)dz\\
=& 1/\lambda_k^2 l_{j}(1)-1/\lambda_k^2 l_{j}(0)\cos(\lambda_k)+
1/\lambda_k^3\int_{0}^1
l'_{j}(z)d\sin\big((1 -z)\lambda_k\big)\\
=& 1/\lambda_k^2 l_{j}(1)-1/\lambda_k^2 l_{j}(0)\cos(\lambda_k)-
1/\lambda_k^3l'_{j}(0)\sin(\lambda_k)-1/\lambda_k^3\int_{0}^1
l''_{j}(z)\sin\big((1 -z)\lambda_k\big)dz\\
=& 1/\lambda_k^2 l_{j}(1)-1/\lambda_k^2 l_{j}(0)\cos(\lambda_k)-
1/\lambda_k^3l'_{j}(0)\sin(\lambda_k)
\\&-1/\lambda_k^4 l''_{j}(1)+1/\lambda_k^4
l''_{j}(0)\cos(\lambda_k)+
1/\lambda_k^5l^{(3)}_{j}(0)\sin(\lambda_k)+1/\lambda_k^5\int_{0}^1
l^{(4)}_{l,j}(z)\sin\big((1 -z)\lambda_k\big)dz\\
=&\cdots\\
=
&\sum\limits_{k=0}^{\lfloor\deg(l_{j})/2\rfloor}(-1)^{k}/\lambda_k^{2k+2}\Big(
l_{j}^{(2k)}(1)- l_{j}^{(2k)}(0)\cos(\lambda_k)-1/\lambda_k
l_{j}^{(2k+1)}(0)\sin(\lambda_k)\Big),\ \ \
i=1,2,\ldots,s,\\
\end{aligned}%
$\\ where $\deg(l_{j})$ is the degree of  $l_{j}$ and
$\lfloor\deg(l_{j})/2\rfloor$ denotes the integral part of
$\deg(l_{j})/2$. Similarly, we obtain
\begin{equation}
\begin{aligned}
&\int_{0}^1l_j(z)\cos\big((1 -z)\lambda_k\big)dz\\=
&\sum\limits_{k=0}^{\lfloor\deg(l_{j})/2\rfloor}(-1)^{k}/\lambda_k^{2k+1}\Big(
l_{j}^{(2k)}(0)\sin(\lambda_k)+1/\lambda_k
l_{j}^{(2k+1)}(1)-1/\lambda_k^2
l_{j}^{(2k+1)}(0)\cos(\lambda_k)\Big),\\
&\int_{0}^1 l_j(c_iz) (c_i\lambda_k )^{-1}\sin\big((1
-z)c_i\lambda_k\big)dz\\=
&\sum\limits_{k=0}^{\lfloor\deg(l_{j})/2\rfloor}(-1)^{k}/(c_i\lambda_k)^{2k+2}\Big(
l_{j}^{(2k)}(c_i)- l_{j}^{(2k)}(0)\cos(c_i\lambda_k)-1/\lambda_k
l_{j}^{(2k+1)}(0)\sin(c_i\lambda_k)\Big),\\
&\qquad i,j=1,2,\ldots,s.\\
\end{aligned}
\label{Q2}%
\end{equation}
\subsection{The scheme of   trigonometric collocation methods} \label{subsec:The methods}
 We are now in a position to present a   kind of trigonometric collocation methods
for the multi-frequency  oscillatory second-order ODEs \eqref{prob}.
\begin{mydef}
\label{numerical method}  A trigonometric collocation method
 for integrating the multi-frequency  oscillatory system \eqref{prob}   is
defined as
\begin{equation}
\begin{aligned} &\tilde{q}_i=\phi_0(c_i^2V)q_0+c_ih\phi_1(c_i^2V)p_0+
(c_ih)^2\sum\limits_{j=1}^ {s}\tilde{I}_{c_i,j}f(\tilde{q}_j ),\quad
i=1,2,\ldots,s,\\
 &\tilde{q}(h)=\phi_0(V)q_0+h\phi_1(V)p_0+ h^2\sum\limits_{j=1}^
{s}I_{1,j}f(\tilde{q}_j),\\
&\tilde{p}(h)=-hM\phi_1( V)q_0+\phi_0(V)p_0 +h\sum\limits_{j=1}^
{s}I_{2,j}f(\tilde{q}_j),
\end{aligned}
\label{methods}%
\end{equation}
where $h$ is the stepsize and  $I_{1,j},\ I_{2,j},\
\tilde{I}_{c_i,j}$ can be computed as stated in Subsection
\ref{subsec:Computations of I}.
\end{mydef}

\begin{rem}
In \cite{wang-2014},  the authors  took advantage of shifted
Legendre polynomials to obtain a local Fourier expansion of the
system \eqref{prob} and derived  trigonometric Fourier collocation
methods (TFCMs). TFCMs are the subclass of $s$-stage ERKN methods
which were presented in  \cite{wu2010-1} with the following Butcher
tableau:
\begin{equation}\label{ERKN tableau-1}
\begin{tabular}
[c]{l}%
\\
\\[2mm]%
$$ $\ $%
\end{tabular}%
\begin{tabular}
[c]{c|ccc}%
$c_{1}$ & $\sum\limits_{j=0}^ {r-1}II_{1,j,c_1}(V)
b_1\widehat{P}_j(c_1)$ & $\ldots$ &
$\sum\limits_{j=0}^ {r-1}II_{1,j,c_1}(V) b_s\widehat{P}_j(c_s)$\\
$\vdots$ & $\vdots$ & $\ddots$ & $\vdots$\\
$c_{s}$ & $\sum\limits_{j=0}^ {r-1}II_{1,j,c_s}(V)
b_1\widehat{P}_j(c_1)$ & $\cdots$ & $\sum\limits_{j=0}^
{r-1}II_{1,j,c_s}(V) b_s\widehat{P}_j(c_s)$\\\hline &
$\raisebox{-1.3ex}[0.5pt]{$\sum\limits_{j=0}^
{r-1}II_{1,j}(V)b_1\widehat{P}_j(c_1)$}$ &
\raisebox{-1.3ex}{$\cdots$} &
$\raisebox{-1.3ex}[0.5pt]{$\sum\limits_{j=0}^
{r-1}II_{1,j}(V)b_s\widehat{P}_j(c_s)$}$\\\hline &
$\raisebox{-1.3ex}[1.0pt]{$\sum\limits_{j=0}^
{r-1}II_{2,j}(V)b_1\widehat{P}_j(c_1)$}$ &
$\raisebox{-1.3ex}[1.0pt]{$\cdots$}$ &
$\raisebox{-1.3ex}[1.0pt]{$\sum\limits_{j=0}^
{r-1}II_{2,j}(V)b_s\widehat{P}_j(c_s)$}$%
\end{tabular}
\end{equation}%
where \begin{equation*}
\begin{aligned} &II_{1,j}(V):=\int_{0}^1\widehat{P}_j(z)(1-z)\phi_1\big((1
-z)^2V\big)dz,\\
&II_{2,j}(V):=\int_{0}^{1}\widehat{P}_j(z)\phi_0\big((1
-z)^2V\big)dz,\\
&II_{1,j,c_i}(V):=\int_{0}^1\widehat{P}_j(c_iz)(1-z)\phi_1\big((1
-z)^2c_i^2V\big)dz,
\end{aligned}
\label{II12j}%
\end{equation*}
  $r$
is an integer with the requirement:  $2\leq r\leq s,$
$\widehat{P}_j$ are shifted Legendre polynomials over the interval
$[0,1]$ and $c_l,\ b_l,\ l=1,2,\ldots,k$ are the node points and the
quadrature weights of a quadrature formula, respectively.

 It is noted that the  method \eqref{methods} is also a subclass of
$s$-stage ERKN methods with the following Butcher tableau:
\begin{equation}\label{tableau-2}%
\begin{tabular}
[c]{l}%
\\
\\[2mm]%
$$ $\ $%
\end{tabular}%
\begin{tabular}
[c]{c|ccc}%
$c_{1}$ & $\tilde{I}_{c_1,1}$ & $\ldots$ &
$\tilde{I}_{c_1,s}$\\
$\vdots$ & $\vdots$ & $\ddots$ & $\vdots$\\
$c_{s}$ & $\tilde{I}_{c_s,1}$ & $\cdots$ &
$\tilde{I}_{c_s,s}$\\\hline & $\raisebox{-1.3ex}[0.5pt]{$I_{1,1}$}$
& \raisebox{-1.3ex}{$\cdots$} &
$\raisebox{-1.3ex}[0.5pt]{$I_{1,s}$}$\\\hline &
$\raisebox{-1.3ex}[1.0pt]{$I_{2,1}$}$ &
$\raisebox{-1.3ex}[1.0pt]{$\cdots$}$ &
$\raisebox{-1.3ex}[1.0pt]{$I_{2,s}$}$%
\end{tabular}
\end{equation}
where
\begin{equation*}
\begin{aligned} &I_{1,j}:=\int_{0}^1l_j(z)(1-z)\phi_1\big((1
-z)^2V\big)dz,\\
&I_{2,j}:=\int_{0}^{1}l_j(z)\phi_0\big((1 -z)^2V\big)dz,\\
&\tilde{I}_{c_i,j}=\int_{0}^1l_j(c_iz)(1-z)\phi_1\big((1
-z)^2c_i^2V\big)dz.
\end{aligned}
\end{equation*}

 From  \eqref{ERKN tableau-1} and
\eqref{tableau-2}, it follows clearly that the coefficients of
\eqref{tableau-2} are simpler than \eqref{ERKN tableau-1}.
Therefore, the scheme of the methods derived in this paper is much
simpler than that given in \cite{wang-2014}. The obtained methods
can be implemented in practical computations at a lower cost, which
will be shown by the numerical experiments in Section \ref{numerical
experiments}. The reason for this point is that we use a more
classical approach and choose Lagrange polynomials to give a local
Fourier expansion of the system \eqref{prob}.
\end{rem}
\begin{rem}
It can be observed from the two tableaus \eqref{ERKN
tableau-1}--\eqref{tableau-2} that the methods here presented are
different from those presented in \cite{wang-2014}. We also note
that in the recent monograph \cite{Brugnano2016}, it has been shown
that the approach of constructing energy-preserving methods for
Hamiltonian problems which are based upon the use of shifted
Legendre polynomials (such as in \cite{Brugnano2012}) and Lagrange
polynomials constructed on Gauss-Legendre nodes (such as in
\cite{hairer2010}) leads to precisely the same methods. Therefore,
by choosing special real numbers $c_1,\cdots,c_s$  for
\eqref{tableau-2} and special quadrature formulae for  \eqref{ERKN
tableau-1}, the methods given in this paper may have some
connections with those in \cite{wang-2014}. We will discuss the
connections in a  future research.
\end{rem}

\begin{rem}
It is noted that the method \eqref{methods} can be applied   to the
system \eqref{prob} with an arbitrary matrix $M$ since the
trigonometric collocation methods do not need the symmetry of $M$.
 Moreover, the method \eqref{methods}  exactly
 integrates the linear system $q''+Mq=0$ and it has an additional advantage of
energy preservation   for linear systems.   The   method
approximates the solution   in the  interval $[0,h]$. We then lend
the procedure with equal ease to  next interval. Namely, we  can
consider the obtained result as the initial condition for a new
initial value problem  in the interval $[h,2h]$. In this way, the
method \eqref{methods} can approximate the solution in an arbitrary
interval $[0,t_{\mathrm{end}}]$  with  $t_{\mathrm{end}}=Nh$.
\end{rem}

When $M= 0$,   \eqref{prob} reduces to  a special and important
class of systems of second-order ODEs  expressed in the traditional
form
\begin{equation}
q^{\prime\prime}(t)=f\big(q(t)\big),  \qquad
q(0)=q_0,\ \ q'(0)=q_0',\qquad t\in[0,t_{\mathrm{end}}].\label{common prob}%
\end{equation}
For this case, with the definition \eqref{Phi01} and the results of
$I_{1,j},\ I_{2,j},\ \tilde{I}_{c_i,j}$   in Subsection
\ref{subsec:Computations of I}, the trigonometric collocation method
\eqref{methods} becomes the following scheme.

\begin{mydef}
\label{numerical RKN method}  An RKN-type    collocation method
 for integrating the traditional second-order ODEs  \eqref{common prob} is
defined as
\begin{equation}
\begin{aligned} &\tilde{q}_i=q_0+c_ihp_0+
(c_ih)^2\sum\limits_{j=1}^
{s}\frac{l_j\Big(\frac{c_i}{3}\Big)}{2}f(\tilde{q}_j),\quad
i=1,2,\ldots,s,\\
 &\tilde{q}(h)=q_0+hp_0+ h^2\sum\limits_{j=1}^
{s}\frac{l_j\Big(\frac{1}{3}\Big)}{2}f(\tilde{q}_j),\\
&\tilde{p}(h)=p_0 +h\sum\limits_{j=1}^
{s}l_j\Big(\frac{1}{2}\Big)f(\tilde{q}_j),
\end{aligned}
\label{methods0}%
\end{equation}
where $h$ is the stepsize.
\end{mydef}

\begin{rem}
 It is noted that the  method \eqref{methods0} is the subclass of
$s$-stage RKN methods with the following Butcher tableau:
\begin{equation}\label{Butcher tableau-RKN}%
\begin{tabular}
[c]{l}%
\\
\\[2mm]%
\begin{tabular}
[c]{c|c}%
$c$ & $\bar{A}=(\bar{a}_{ij})_{k\times k}$\\\hline &
$\raisebox{-1.3ex}[0pt]{$\bar{b}^T$}$\\\hline
& $\raisebox{-1.3ex}[0.5pt]{$b^T$}$%
\end{tabular}
$\ \quad=$ $\ $%
\end{tabular}%
\begin{tabular}
[c]{c|ccc}%
$c_{1}$ & $l_1\Big(\frac{c_1}{3}\Big)/2$ & $\ldots$ &
$l_s\Big(\frac{c_1}{3}\Big)/2$\\
$\vdots$ & $\vdots$ & $\ddots$ & $\vdots$\\
$c_{s}$ & $l_1\Big(\frac{c_s}{3}\Big)/2$ & $\cdots$ &
$l_s\Big(\frac{c_s}{3}\Big)/2$\\\hline &
$\raisebox{-1.3ex}[0.5pt]{$l_1\Big(\frac{1}{3}\Big)/2$}$ &
\raisebox{-1.3ex}{$\cdots$} &
$\raisebox{-1.3ex}[0.5pt]{$l_s\Big(\frac{1}{3}\Big)/2$}$\\\hline &
$\raisebox{-1.3ex}[1.0pt]{$l_1\Big(\frac{1}{2}\Big)$}$ &
$\raisebox{-1.3ex}[1.0pt]{$\cdots$}$ &
$\raisebox{-1.3ex}[1.0pt]{$l_s\Big(\frac{1}{2}\Big)$}$%
\end{tabular}
\end{equation}
 This point means that by letting $M=0$, the   trigonometric collocation methods yield a subclass of
  RKN methods for solving traditional second-order ODEs,
which  demonstrates the wider applications of the  methods.
\end{rem}

\section{Properties of the    methods} \label{sec:Analysis of the methods}

For the exact solution of \eqref{H-s} at $t=h$,  let
$\mathbf{y}(h)=\Big(
                        q^{\intercal}(h),p^{\intercal}(h)\Big)^{\intercal}.$
 Then the oscillatory Hamiltonian system \eqref{H-s} can be rewritten
                   in the form
\begin{equation}\mathbf{y}'(\xi h)=F(\mathbf{y}(\xi h)):=\left(
                                                                           \begin{array}{c}
                                                                             p(\xi h) \\
                                                                            -Mq(\xi
                                                                            h)+f\big(q(\xi h)\big)
                                                                           \end{array}
                                                                         \right),\quad
\mathbf{y}_0=\left(
                 \begin{array}{c}
                   q_0 \\
                   p_0 \\
                 \end{array}
               \right)
\label{new prob}%
\end{equation}
for $0\leq\xi \leq1.$
 The Hamiltonian is \begin{equation}H(\mathbf{y})=\frac{1}{2}p^{\intercal}p+\frac{1}{2}q^{\intercal}Mq+U(q).
 \label{new H}%
\end{equation}

  On the other hand,  denoting the numerical method
\eqref{methods} as $$\mathbf{\omega}(h)=\Big(
                        \tilde{q}^{\intercal}(h),
                     \tilde{p}^{\intercal}(h)\Big)^{\intercal},$$  the numerical solution
                   satisfies
\begin{equation}\mathbf{\omega}'(\xi h)=\left(
                                                                           \begin{array}{c}
                                                                             \tilde{p}(\xi h) \\
                                                                           -M\tilde{q}(\xi h)+\sum\limits_{j=1}^
{s}l_j(\xi )f\big(\tilde{q}(c_j h)\big)
                                                                           \end{array}
                                                                         \right),\quad
\mathbf{\omega}_0=\left(
                 \begin{array}{c}
                   q_0 \\
                   p_0 \\
                 \end{array}
               \right).
\label{new solver}%
\end{equation}
 The next lemma is useful for the following analysis.
\begin{lem}
\label{Pj lem}Let $g:[0,h]\rightarrow \mathbb{R}^{d}$ have $j$
continuous derivatives in the interval $[0,h]$. Then
$$\int_{0}^1P_j(\tau)g(\tau h)d\tau=\mathcal{O}(h^{j}),$$ where
$P_j(\tau)$ is an orthogonal polynomial of degree $j$ on the
interval $[0,1]$.
\end{lem}
\textbf{Proof.}
  We assume that  $g(\tau h)$ can be expanded in Taylor
series at the origin for sake of simplicity.  Then, for all
$j\geq0$, by considering that $P_j(\tau)$ is orthogonal to all
polynomials of degree $n< j$:
$$\int_{0}^1P_j(\tau)g(\tau h)d\tau=\sum\limits_{n=1}^
{\infty}\frac{g^{(n)}(0)}{n!}h^n\int_{0}^1P_j(\tau)\tau^nd\tau=\mathcal{O}(h^{j}).$$
\hfill  $\blacksquare$
%
\subsection{The order of  energy preservation} \label{subsec:The second invariant preserving of the methods}
In this subsection we are concerned with the  order of  preservation
of   the Hamiltonian energy.
\begin{mytheo}\label{energy thm}
Under  the condition that  $c_l,\ l=1,2,\ldots,s$ are chosen as the
node points of  a $s$-point Gauss--Legendre's
 quadrature over the integral $[0,1]$,   we have
$$H(\omega(h))-H(\mathbf{y}_0)=\mathcal{O}(h^{2s+1}),$$
 where the constant symbolized
by   $\mathcal{O}$    is independent of   $h$.
\end{mytheo}
\textbf{Proof.}  By virtue of Lemma \ref{Pj lem}, \eqref{new H} and
\eqref{new solver}, one has
\begin{equation*}
\begin{aligned} &H(\omega(h))-H(\mathbf{y}_0)
=h \int_{0}^{1}
\nabla H(\omega(\xi  h))^{\intercal}\omega'(\xi  h)d\xi \\
=&h \int_{0}^{1} \Big( \big(M\tilde{q}(\xi h)-f(\tilde{q}(\xi
h)\big)^{\intercal},\ \tilde{p}(\xi
h)^{\intercal}\Big)\\
& \cdot\left(
   \begin{array}{c}
  \tilde{p}(\xi h) \\
   -M\tilde{q}(\xi h)+\sum\limits_{j=1}^ {s}l_j(\xi
)f\big(\tilde{q}(c_j h)\big)
 \end{array}
  \right)d\xi \\
=&h \int_{0}^{1} \tilde{p}(\xi h)^{\intercal} \Big(
\sum\limits_{j=1}^ {s}l_j(\xi )f(\tilde{q}(c_j
h))-f\big(\tilde{q}(\xi h)\big) \Big)d\xi.
\end{aligned}
\end{equation*}
Moreover, we have
$$f\big(\tilde{q}(\xi h)\big)-\sum\limits_{j=1}^ {s}l_j(\xi )f\big(\tilde{q}(c_j
h)\big)=\frac{f^{(s+1)}\big(\tilde{q}(\xi
h)\big)|_{\xi=\zeta}}{(n+1)!}\prod\limits_{i=1}^ {s}(\xi h-c_ih).$$
Here $f^{(s+1)}\big(\tilde{q}(\xi h)\big)$ denote the
$(s+1)$th-order derivative of $f(\tilde{q}(t))$ with respect to $t$.
Then,  we obtain
\begin{equation*}
\begin{aligned} H(\omega(h))-H(\mathbf{y}_0)
=&-h \int_{0}^{1} \tilde{p}(\xi h)^{\intercal}
\frac{f^{(s+1)}\big(\tilde{q}(\xi
h)\big)|_{\xi=\zeta}}{(n+1)!}\prod\limits_{i=1}^
{s}(\xi h-c_ih)  d\xi \\
=&-h^{s+1}\int_{0}^{1} \tilde{p}(\xi h)^{\intercal}
\frac{f^{(s+1)}\big(\tilde{q}(\xi
h)\big)|_{\xi=\zeta}}{(n+1)!}\prod\limits_{i=1}^ {s}(\xi -c_i) d\xi.
\end{aligned}
\end{equation*}

Since $c_l,\ l=1,2,\ldots,s$ are chosen as the node points of  a
$s$-point Gauss--Legendre's
 quadrature over the integral $[0,1]$, $\prod\limits_{i=1}^ {s}(\xi
-c_i)$ is an orthogonal polynomial of degree $s$ on the interval
$[0,1]$. Therefore, it follows from Lemma \ref{Pj lem} that
\begin{equation*}
\begin{aligned} &H(\omega(h))-H(\mathbf{y}_0)
=h^{s+1}\mathcal{O}(h^{s})=\mathcal{O}(h^{2s+1}).\\
\end{aligned}
\end{equation*}
\hfill  $\blacksquare$

\subsection{The   order of  quadratic invariant} \label{subsec:The second invariant preserving of the methods}
We  next turn to  the quadratic invariant
$Q(\mathbf{y})=q^{\intercal}Dp$ of \eqref{prob}. The quadratic form
$Q$ is a first integral of \eqref{prob} if and only if
$p^{\intercal}Dp+q^{\intercal}D(f(q)-Mq)=0$   for all
$p,q\in\mathbb{R}^{d}$. This implies that   $D$ is a skew-symmetric
matrix and that $q^{\intercal}D(f(q)-Mq)=0$ for any
$q\in\mathbb{R}^{d}$. The following result states the degree of
accuracy of the  method \eqref{methods}.
\begin{mytheo}\label{invariant thm}
Under the condition in Theorem \ref{energy thm}, we have
$$Q(\omega(h))-Q(\mathbf{y}_0)=\mathcal{O}(h^{2s+1}),$$
 where the constant symbolized
by   $\mathcal{O}$    is independent of  $h$.
\end{mytheo}
 \textbf{Proof.} From  $Q(\mathbf{y})=q^{\intercal}Dp$ and
$D^{\intercal}=-D$, it follows that
\begin{equation*}
\begin{aligned} &Q(\omega(h))-Q(\mathbf{y}_0)
=h \int_{0}^{1}
\nabla Q(\omega(\xi  h))^{\intercal}\omega'(\xi  h)d\xi \\
=&h \int_{0}^{1} \Big(- \tilde{p}(\xi h)^{\intercal}D,\
\tilde{q}(\xi h)^{\intercal}D\Big)\left(
   \begin{array}{c}
  \tilde{p}(\xi h) \\
   -M\tilde{q}(\xi h)+\sum\limits_{j=1}^ {s}l_j(\xi
)f\big(\tilde{q}(c_j h)\big)
 \end{array}
  \right)d\xi .
\end{aligned}
\end{equation*}
Since   $q^{\intercal}D(f(q)-Mq)=0$ for any $q\in\mathbb{R}^{d}$, we
obtain
\begin{equation*}
\begin{aligned} &Q(\omega(h))-Q(\mathbf{y}_0)
  =h \int_{0}^{1} \tilde{q}(\xi h)^{\intercal}D
   \Big(-M\tilde{q}(\xi h)+\sum\limits_{j=1}^ {s}l_j(\xi
)f\big(\tilde{q}(c_j h)\big)\Big)d\xi \\
     =&h \int_{0}^{1} \tilde{q}(\xi h)^{\intercal}D
   \frac{f^{(s+1)}\big(\tilde{q}(\xi
h)\big)|_{\xi=\zeta}}{(n+1)!}\prod\limits_{i=1}^
{s}(\xi h-c_ih)d\xi \\
     =&h^{s+1} \int_{0}^{1} \tilde{q}(\xi h)^{\intercal}D
   \frac{f^{(s+1)}\big(\tilde{q}(\xi
h)\big)|_{\xi=\zeta}}{(n+1)!}\prod\limits_{i=1}^
{s}(\xi -c_i)d\xi \\
 =&\mathcal{O}(h^{s+1})\mathcal{O}(h^{s})=\mathcal{O}(h^{2s+1}).
\end{aligned}
\end{equation*}\hfill  $\blacksquare$

\subsection{The order} \label{subsec:The order of the methods}
To express the dependence of the solutions of
$\mathbf{y}'(t)=F(\mathbf{y}(t))$ on the initial values, for any
given $\tilde{t}\in[0,h]$, we  denote by
$\mathbf{y}(\cdot,\tilde{t}, \tilde{\mathbf{y}})$ the solution
satisfying the initial condition $\mathbf{y}(\tilde{t},\tilde{t},
\tilde{\mathbf{y}})=\tilde{\mathbf{y}}$ and   set
\begin{equation}
\Phi(s,\tilde{t},
\tilde{\mathbf{y}})=\frac{\partial\mathbf{y}(s,\tilde{t},
\tilde{\mathbf{y}})}{\partial \tilde{\mathbf{y}}}.
\label{Phi}%
\end{equation}
Recalling the elementary theory of ODEs, we have the following
standard result   (see, e.g.  \cite{hale1980})
\begin{equation}
\frac{\partial\mathbf{y}(s,\tilde{t}, \tilde{\mathbf{y}})}{\partial
\tilde{t}}=-\Phi(s,\tilde{t},
\tilde{\mathbf{y}})F(\tilde{\mathbf{y}}).
\label{standard result}%
\end{equation}

The following theorem states the result on the  order of the
trigonometric collocation methods.
\begin{mytheo}\label{order thm}
Under the condition in Theorem \ref{energy thm},  the trigonometric
collocation
 method \eqref{methods} satisfies
$$\mathbf{y}(h)-\omega(h)=\mathcal{O}(h^{2s+1}),$$
 where the constant symbolized
by   $\mathcal{O}$    is independent of   $h$.

\end{mytheo}
\textbf{Proof.}  It follows from   \eqref{Phi} and \eqref{standard
result} that
\begin{equation*}
\begin{aligned} &\mathbf{y}(h)-\omega(h)
=\mathbf{y}(h,0, \mathbf{y}_0)-\mathbf{y}\big(h,h, \omega(h)\big)=-
\int_{0}^{h}
\frac{d\mathbf{y}\big(h,\tau, \omega(\tau)\big)}{d\tau}d\tau\\
=& - \int_{0}^{h}\Big[ \frac{\partial\mathbf{y}\big(h,\tau,
\omega(\tau)\big)}{\partial \tilde{t}}
+\frac{\partial\mathbf{y}\big(h,\tau, \omega(\tau)\big)}{\partial \tilde{\mathbf{y}}}\omega'(\tau)\Big]d\tau\\
=& h \int_{0}^{1}\Phi\big(h,\xi h, \omega(\xi h)\big)\Big[F\big(\omega(\xi h)\big)-\omega'(\xi h)\Big]d\xi \\
=& h \int_{0}^{1}\Phi\big(h,\xi h, \omega(\xi h)\big) \left(
                                                                           \begin{array}{c}
                                                                             \mathbf{0}\\
f\big(\tilde{q}(\xi h)\big)-\sum\limits_{j=1}^ {s}l_j(\xi
)f\big(\tilde{q}(c_j h)\big)
                                                                           \end{array}
                                                                         \right)d\xi .
\end{aligned}
\end{equation*}
We rewrite $\Phi\big(h,\xi h, \omega(\xi h)\big)$ as a block
 matrix:
$$\Phi\big(h,\xi h, \omega(\xi h)\big)=\left(
                           \begin{array}{cc}
                             \Phi_{11}(\xi h) & \Phi_{12}(\xi h) \\
                             \Phi_{21}(\xi h) & \Phi_{22}(\xi h) \\
                           \end{array}
                         \right),
$$
where $\Phi_{ij}\ (i,j=1,2)$ are all $d\times d$ matrices.

 We then yield
\begin{equation*}
\begin{aligned} &\mathbf{y}(h)-\omega(h)
=h  \left(
                                                                           \begin{array}{c}
\int_{0}^{1}\Phi_{12}(\xi h)\frac{f^{(s+1)}\big(\tilde{q}(\xi
h)\big)|_{\xi=\zeta}}{(n+1)!}\prod\limits_{i=1}^ {s}(\xi
h-c_ih)d\xi \\
                                                                            \int_{0}^{1}\Phi_{22}(\xi h)
                                                                            \frac{f^{(s+1)}\big(\tilde{q}(\xi
h)\big)|_{\xi=\zeta}}{(n+1)!}\prod\limits_{i=1}^ {s}(\xi h-c_ih)d\xi
                                                                           \end{array}
                                                                         \right)\\
=&h^{s+1}  \left(
                                                                           \begin{array}{c}
\int_{0}^{1}\Phi_{12}(\xi h)\frac{f^{(s+1)}\big(\tilde{q}(\xi
h)\big)|_{\xi=\zeta}}{(n+1)!}\prod\limits_{i=1}^ {s}(\xi
-c_i)d\xi \\
                                                                            \int_{0}^{1}\Phi_{22}(\xi h)\frac{f^{(s+1)}\big(\tilde{q}(\xi
h)\big)|_{\xi=\zeta}}{(n+1)!}\prod\limits_{i=1}^ {s}(\xi -c_i)d\xi
                                                                           \end{array}
                                                                         \right)=h^{s+1}\mathcal{O}(h^{s})      =\mathcal{O}(h^{2s+1}).\\
\end{aligned}
\end{equation*}
\hfill  $\blacksquare$

\subsection{Convergence analysis of the iteration} \label{subsec:Convergence}

\begin{mytheo}
\label{convergence} Assume that $M$ is symmetric and positive
semi-definite and that  $f$ satisfies a Lipschitz condition in the
variable $q$, i.e., there exists a constant $L$ with the property
that $\norm{f(q_1)-f(q_2)}\leq L\norm{q_1-q_2}$. If
\begin{equation}0<h<\frac{1}{\sqrt{L\max\limits_{i,j= 1,\cdots, s}\int_{0}^1|l_j(c_iz)(1-z)|dz}}, \label{Convergence
condition}
\end{equation} then the fixed-point  iteration for
the  method   \eqref{methods} is convergent.
\end{mytheo}
\textbf{Proof.}  Following Definition \ref{numerical method}, the
first formula
  of   \eqref{methods} can be rewritten as
\begin{align}
Q  & =\phi_{0}(c^{2}V)q_{0}+c\phi_{1}(c^{2}V)hp_{0}+h^2A(V)f(Q),
\end{align}
where $c=(c_1,c_2,\ldots,c_k)^{\intercal},\
Q=(\tilde{q}_1,\tilde{q}_2,\ldots,\tilde{q}_k)^{\intercal},\
A(V)=\big(a_{ij}(V)\big)_{k\times k}$ and $a_{ij}(V)$ are defined as
$$a_{ij}(V):= \int_{0}^1l_j(c_iz)(1-z)\phi_1\big((1
-z)^2c_i^2V\big)dz.
$$

By Proposition 2.1 in \cite{2012_li}, we know that
$\norm{\phi_1\big((1 -z)^2c_i^2V\big)}\leq1$ and then we get
 \begin{equation*}
\begin{aligned}\norm{a_{ij}(V)}
&\leq  \int_{0}^1|l_j(c_iz)(1-z)|dz.\\
\end{aligned}
\end{equation*}

 Let
$$\varphi(x)=\phi_{0}(c^{2}V)q_{0}+c\phi_{1}(c^{2}V)hp_{0}+h^2A(V)f(x).$$ Then
\begin{equation*}
\begin{aligned}
\norm{\varphi(x)-\varphi(y)} &=\norm{h^2A(V)f(x)-h^2A(V)f(y)} \leq
 h^2L\norm{A(V)} \norm{x-y}\\
 &\leq h^2L\max\limits_{i,j= 1,\cdots, s}\int_{0}^1|l_j(c_iz)(1-z)|dz\norm{x-y},
\end{aligned}
\end{equation*}
which means that $\varphi(x)$ is a contraction from the assumption
\eqref{Convergence condition}. The well-known Contraction Mapping
Theorem then ensures the convergence of the fixed-point iteration.
$\blacksquare$

\begin{rem}
It is noted that the   convergence of the   methods is independent
of $\norm{M}$. This point is  of prime importance especially for
highly oscillatory systems since  we usually have $\norm{M}\gg1$,
which will be  shown by the numerical results of Problem 2 in
Section \ref{numerical experiments}.
\end{rem}

\subsection{Stability and phase
properties} \label{subsec:stability} In this part we are concerned
with the stability and phase properties.　We consider the test
equation:
\begin{equation}
q^{\prime\prime}(t)+\omega^{2}q(t)=-\epsilon q(t)\ \ \mathrm{with}
\ \ \omega^{2}+ \epsilon>0, \label{ts}%
\end{equation}
where $\omega$ represents an estimation of the dominant frequency
$\lambda$ and $\epsilon=\lambda^{2}-\omega^{2}$ is the error of that
estimation. Applying \eqref{methods} to \eqref{ts} produces

\[
\left(
\begin{array}
[c]{c}%
\tilde{q}\\
h\tilde{p}%
\end{array}
\right)  =S(V,z)\left(
\begin{array}
[c]{c}%
q_{0}\\
hp_{0}%
\end{array}
\right)  ,
\]
where the stability matrix $S(V,z)$ is given by
\[
S(V,z)=\left(
\begin{array}
[c]{cc}%
\phi_{0}(V)-z\bar{b}^{\intercal}(V)N^{-1}\phi_{0}(c^{2}V) & \phi_{1}(V)\!-\!z\bar{b}%
^{\intercal}(V)N^{-1}(c\cdot\phi_{1}(c^{2}V))\\
-V\phi_{1}(V)\!-\!zb^{\intercal}(V)N^{-1}\phi_{0}(c^{2}V) & \phi_{0}(V)\!-\!zb^{\intercal}%
(V)N^{-1}(c\cdot\phi_{1}(c^{2}V))
\end{array}
\right)
\]
with $N=I+zA(V)$,
  $\bar{b}(V)=\Big(I_{1,1},\ldots,I_{1,s} \Big)^{\intercal},\   b(V)=\Big(I_{2,1},\ldots,I_{2,s}
\Big)^{\intercal}.$

Accordingly, we have the following definitions of stability and
dispersion order and dissipation order for our method
\eqref{methods}.

\begin{mydef}
\label{stability region} (\cite{wu2012-2amm})Let $\rho(S)$ be the
spectral radius of $S$. $$R_{s}=\{(V,z)|\ V>0\ \textmd{and}\
\rho(S)<1\}$$ is called the \textit{stability region of the method
\eqref{methods}}. $$R_{p}=\{(V,z)|\ V>0,\ \rho(S)=1\ \textmd{and}\
\mathrm{tr}(S)^{2}<4\det(S)\}$$ is called the \textit{periodicity
region of the method \eqref{methods}}. The quantities
$$\phi(\zeta)=\zeta-\arccos\Big(\frac
{\mathrm{tr}(S)}{2\sqrt{\det(S)}}\Big),\ \
d(\zeta)=1-\sqrt{\det(S)}$$ are called the dispersion error and the
dissipation error of the method \eqref{methods}, respectively, where
$\zeta=\sqrt{V+z}$. Then, a method is said to be dispersive of order
$r$ and dissipative of order $s$, if $\phi(\zeta)=O(\zeta^{r+1})$
and $d(\zeta)=O(\zeta^{s+1})$, respectively. If $\phi(\zeta)=0$ and
$d(\zeta)=0$, then the method is said to be zero dispersive and zero
dissipative, respectively.
\end{mydef}

 \section{Numerical experiments}
 \label{numerical experiments} As an example of the   trigonometric collocation methods \eqref{methods}, we   choose   the node points of  a two-point
Gauss--Legendre's quadrature over the integral $[0,1]$
\begin{equation}\begin{aligned}
&c_1=\frac{3-\sqrt{3}}{6},\ \ c_2=\frac{3+\sqrt{3}}{6}.  \label{GL}
\end{aligned}\end{equation} Then we choose $s=2$
in \eqref{methods} and denote the corresponding  fourth-order method
as LTCM.

\begin{figure}[t!]
\centering%
\begin{tabular}
[c]{ll}%
\includegraphics[width=6cm,height=5cm]
{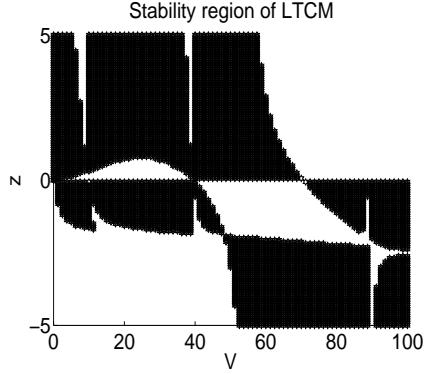} &
\end{tabular}
\caption{{\protect\small Stability region (shaded area) of the method LTCM.}}%
\label{sta}%
\end{figure}

The stability region  of this method is shown in Fig. \ref{sta}. We
note that in order to obtain any information for the stability
regions, we need to consider various values of $V$ and $z$.  Here we
choose the subsets $V\in[0,100],\ z\in[-5,5]$ and these regions
shown in Figure \ref{sta} only give an indication of the stability
of this method.

The dissipative error and dispersion error  are given respectively
by
\begin{equation*}\begin{aligned}
d(\zeta) &= \frac{\epsilon^2}{24(\epsilon+\omega
^{2})^2}\zeta^4+\mathcal{O}(\zeta^{5}),\ \ \ \ \phi(\zeta) =
\frac{\epsilon^2}{6(\epsilon+\omega
^{2})^2}\zeta^3+\mathcal{O}(\zeta^{4}).
\end{aligned}\end{equation*}
It is noted that when   $M=0$,   the   method LTCM reduces to a
fourth-order RKN method with the   Butcher tableau \eqref{Butcher
tableau-RKN} and \eqref{GL}.

In order to   show the efficiency and robustness of the fourth-order
method, the other integrators we select for comparison are:

\begin{itemize}\itemsep=-0.2mm

\item TFCM: a fourth-order trigonometric  Fourier collocation method in
\cite{wang-2014} with $c_1=\frac{3-\sqrt{3}}{6},\
c_2=\frac{3+\sqrt{3}}{6},\ b_1=b_2=1/2,\ r=2$;

\item SRKM1: the symplectic Runge--Kutta method of order five  in \cite{Sun1993} based on Radau
quadrature;

\item EPCM1:  the ``extended Lobatto IIIA
method of order four" in \cite{Iavernaro2009}, which is an
energy-preserving collocation method (the case $s=2$ in
\cite{hairer2010});

\item EPRKM1: the energy-preserving Runge--Kutta method of order four (formula (19) in
\cite{Brugnano2012}).

\end{itemize}

Since all these methods are implicit, we use the classical waveform
Picard algorithm.  For each experiment, first we  show the
convergence rate of iterations  for different error tolerances. Then
for different methods,  we set the error tolerance as $10^{-16}$ and
set the maximum number of iteration as 5. We display the global
errors and the energy errors if the problem is a Hamiltonian system.
The numerical experiments have been carried out on a personal
computer and the algorithm has been implemented by using the
MATLAB-R2013a.

\begin{figure}[ptbh]
\centering\tabcolsep=4mm
\begin{tabular}
[c]{cc}%
\includegraphics[width=5cm,height=6cm]{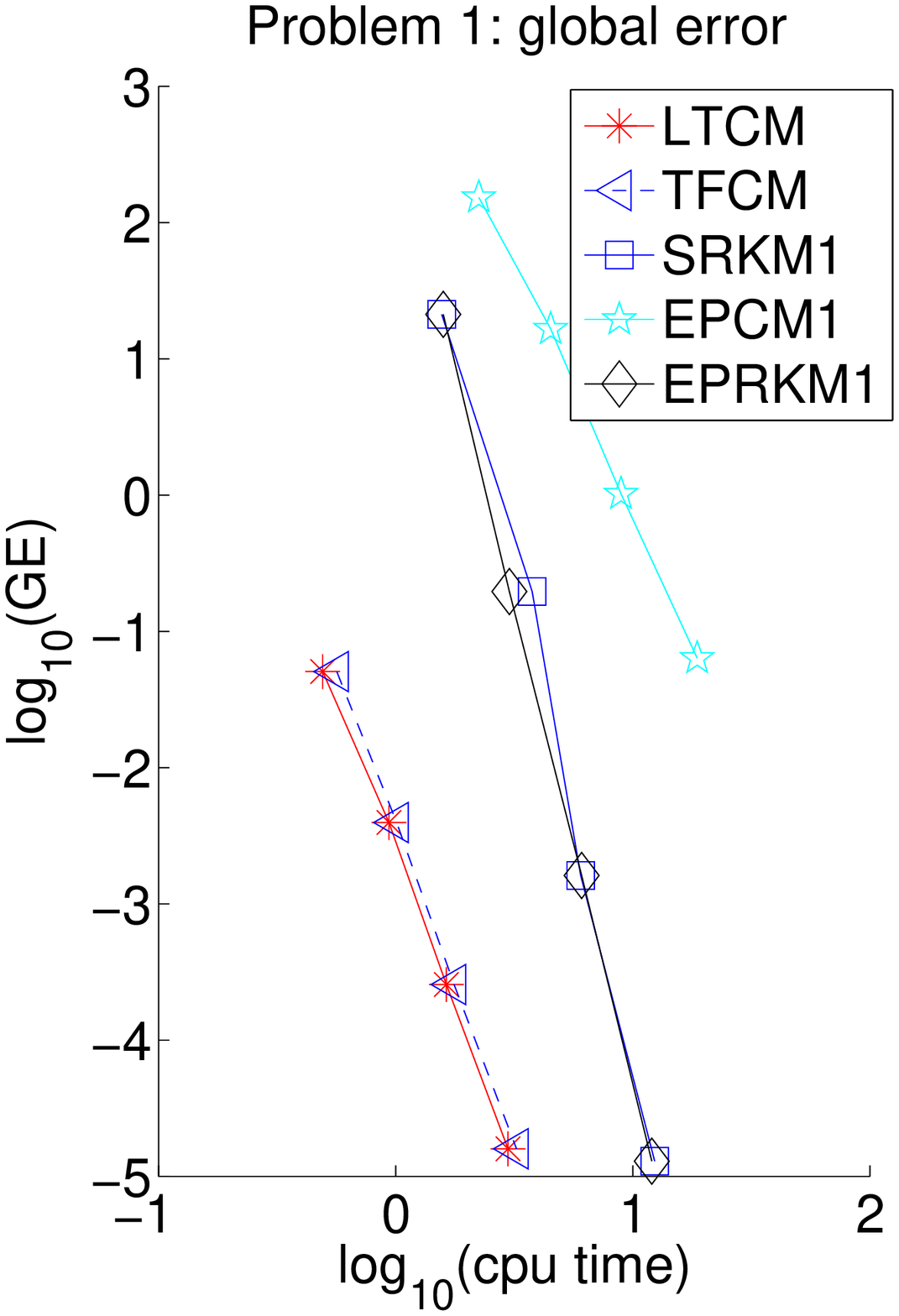} & \includegraphics[width=5cm,height=6cm]{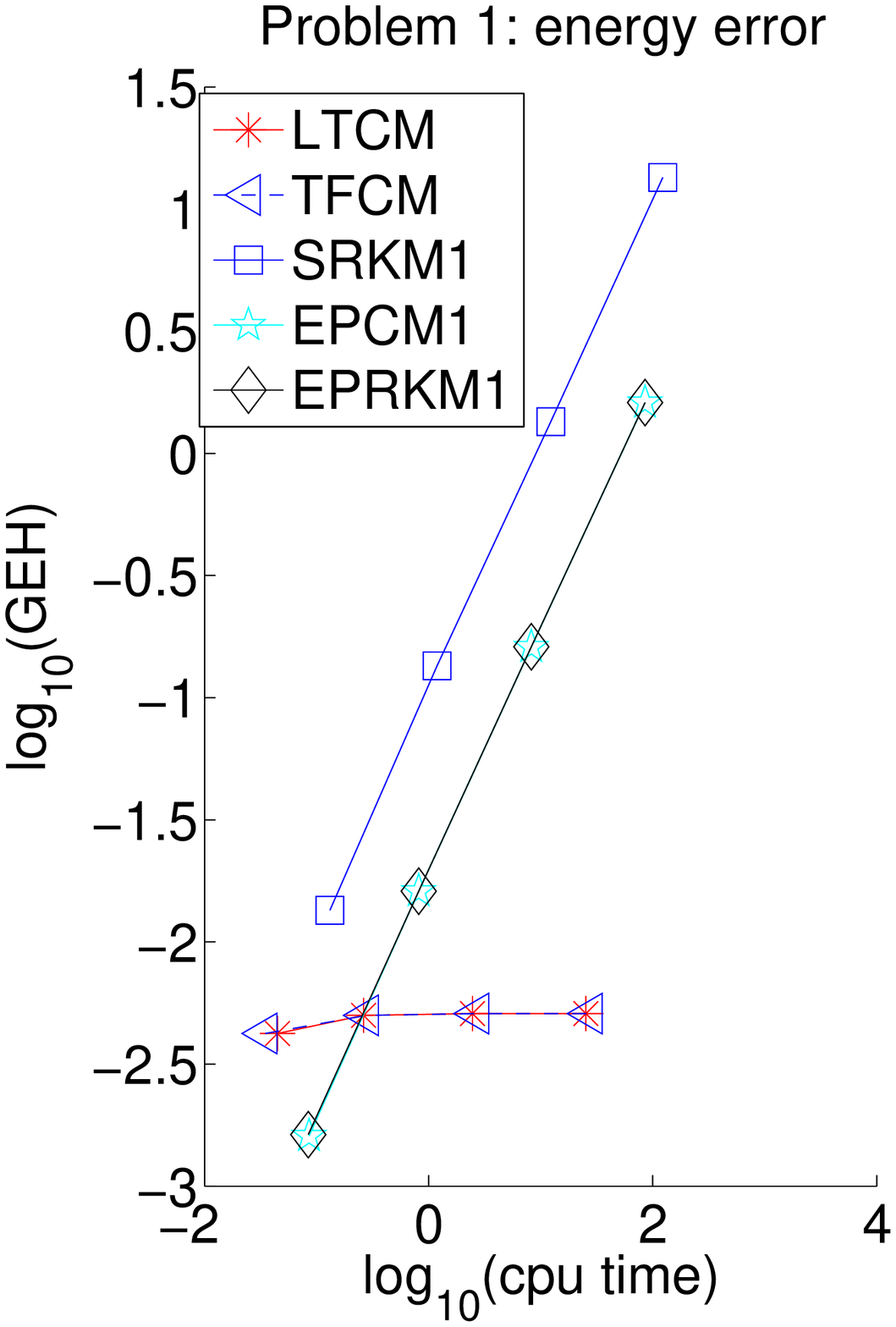}\\
{\small (i)} & {\small (ii)}%
\end{tabular}
\caption{Results for Problem 1. (i): The logarithm of the global
error ($GE$) over the integration interval against the logarithm of
CPU time. (ii):\ The logarithm of the   maximum global error of
Hamiltonian energy ($GEH$)  against the logarithm of
CPU time.}%
\label{fig:problem1}%
\end{figure}

\vskip0.3cm\noindent\textbf{Problem 1.} Consider the Hamiltonian
equation which governs the motion of an artificial satellite  (this
problem has been considered in \cite{stiefel}) with the Hamiltonian

\[H(q,p)=\frac{1}{2}p^{\intercal}p+\frac{1}{2}\frac{\kappa}{2}q^{\intercal}q+\lambda
\Big(\frac{(q_1q_3+q_2q_4)^2}{r^4}-\frac{1}{12r^2}\Big),\] where
$q=(q_1,q_2,q_3,q_4)^{\intercal}$  and $r=q^{\intercal}q.$  The
initial conditions are given on an elliptic equatorial orbit by
$$q_0=\sqrt{\frac{r_0}{2}}\Big(-1,-\frac{\sqrt{3}}{2},-\frac{1}{2},0\Big)^{\intercal},\
\ \
p_0=\frac{1}{2}\sqrt{K^2\frac{1+e}{2}}\Big(1,\frac{\sqrt{3}}{2},\frac{1}{2},0\Big)^{\intercal}.$$
Here $M=\frac{\kappa}{2}$ and $\kappa$ is the total energy of the
elliptic motion which is defined by
$\kappa=\frac{K^2-2|p_0|^2}{r_0}-V_0 $ with
$V_0=-\frac{\lambda}{12r_0^3}.$ The parameters of this problem are
chosen as
 $K^2=3.98601\times10^5$,
$r_0=6.8\times10^3$, $e=0.1$,  $\lambda=\frac{3}{2}K^2J_2R^2,\
J_2=1.08625\times10^{-3},\ R=6.37122\times10^3$. First the problem
is solved in the interval $[0, 10^4]$ with the stepsize
$h=\frac{1}{10}$ to show the convergence rate of iterations. See
Table \ref{pro1-NEW tab} for the CPU time of iterations  for
different error tolerances. Then this equation is integrated in
$[0,1000]$   with the stepsizes $1/2^i$, $i=2,\cdots,5$. The global
errors against CPU time are shown in Fig.  \ref{fig:problem1} (i).
We finally integrate this problem with a fixed stepsize  $h=1/20$ in
the interval $[0,t_{\mathrm{end}}]$ with $t_{\mathrm{end}}=10, 100,
10^3, 10^4$. The maximum global errors of Hamiltonian energy against
CPU time are presented in Fig.  \ref{fig:problem1} (ii).

\begin{table}$$
\begin{array}{|c|c|c|c|c|c|}
\hline
\text{Methods} &tol=1.0e-006  &tol=1.0e-008   &tol=1.0e-010 &tol=1.0e-012    \\
\hline
\text{LTCM} & 6.8215 &   8.8964  &  8.8500  & 10.5551 \cr
\text{TFCM} & 9.7892  &  9.7553 &   9.9806  & 13.0105 \cr
\text{SRKM1} &67.0230  & 64.1777 &  75.9390 &  86.8317\cr
\text{EPCM1}  &104.4341 & 112.9710 & 126.4438  &145.6188\cr
\text{EPRKM1}  & 56.2409 &  64.3123&   75.2503 &  84.9962\cr
 \hline
\end{array}
$$
\caption{Results for Problem 1: The total CPU time (s) of iterations
for different error tolerances (tol).} \label{pro1-NEW tab}
\end{table}

\vspace{2mm} \noindent\textbf{Problem 2.} Consider the
Fermi-Pasta-Ulam Problem \cite{hairer2006}.

Fermi-Pasta-Ulam Problem is a Hamiltonian system with the
Hamiltonian
\[%
\begin{array}
[c]{ll}%
H(y,x) & =\frac{1}{2}\textstyle\sum\limits_{i=1}^{2m}y_{i}^{2}+\frac
{\omega^{2}}{2}\textstyle\sum\limits_{i=1}^{m}x_{m+i}^{2}+\frac{1}{4}%
\Big[(x_{1}-x_{m+1})^{4}\\
& +\textstyle\sum\limits_{i=1}^{m-1}(x_{i+1}-x_{m+i-1}-x_{i}-x_{m+i}%
)^{4}+(x_{m}+x_{2m})^{4}\Big],
\end{array}
\]
where $x_{i}$ is a scaled displacement of the $i$th stiff spring,
$x_{m+i}$ represents a scaled expansion (or compression) of the
$i$th stiff spring, and $y_{i},\ y_{m+i}$ are their velocities (or
momenta). This system can be rewritten as
\[
x^{\prime\prime}(t)+Mx(t)=-\nabla U(x),\qquad t\in\lbrack
t_{0},t_{\mathrm{end}}],
\]
where
\begin{gather*}
M=\left(
\begin{array}
[c]{cc}%
\mathbf{0}_{m\times m} & \mathbf{0}_{m\times m}\\
\mathbf{0}_{m\times m} & \omega^{2}I_{m\times m}%
\end{array}
\right),\\
U(x)=\frac{1}{4}\Big[(x_{1}-x_{m+1})^{4}+\textstyle\sum\limits_{i=1}%
^{m-1}(x_{i+1}-x_{m+i-1}-x_{i}-x_{m+i})^{4}+(x_{m}+x_{2m})^{4}\Big].
\end{gather*}

Following \cite{hairer2006}, we choose
\[
m=3,\ x_{1}(0)=1,\ y_{1}(0)=1,\ x_{4}(0)=\frac{1}{\omega},\
y_{4}(0)=1
\]
with zero for the remaining initial values.

First the problem is solved in the interval $[0, 1000]$ with the
stepsize $h=\frac{1}{100}$ and $\omega=100,\ 200$ to show the
convergence rate of iterations. See Table \ref{pro2-NEW tab-1}  for
the total
 CPU time of iterations for different error tolerances. It can be observed that when $\omega$
increases, the convergence rate of LTCM and TFCM is almost
unaffected. However, the convergence rate of the other methods
varies greatly when $\omega$ becomes large.

\begin{figure}[ptbh]
\centering\tabcolsep=4mm
\begin{tabular}
[c]{cccc}%
\includegraphics[width=5cm,height=6cm]{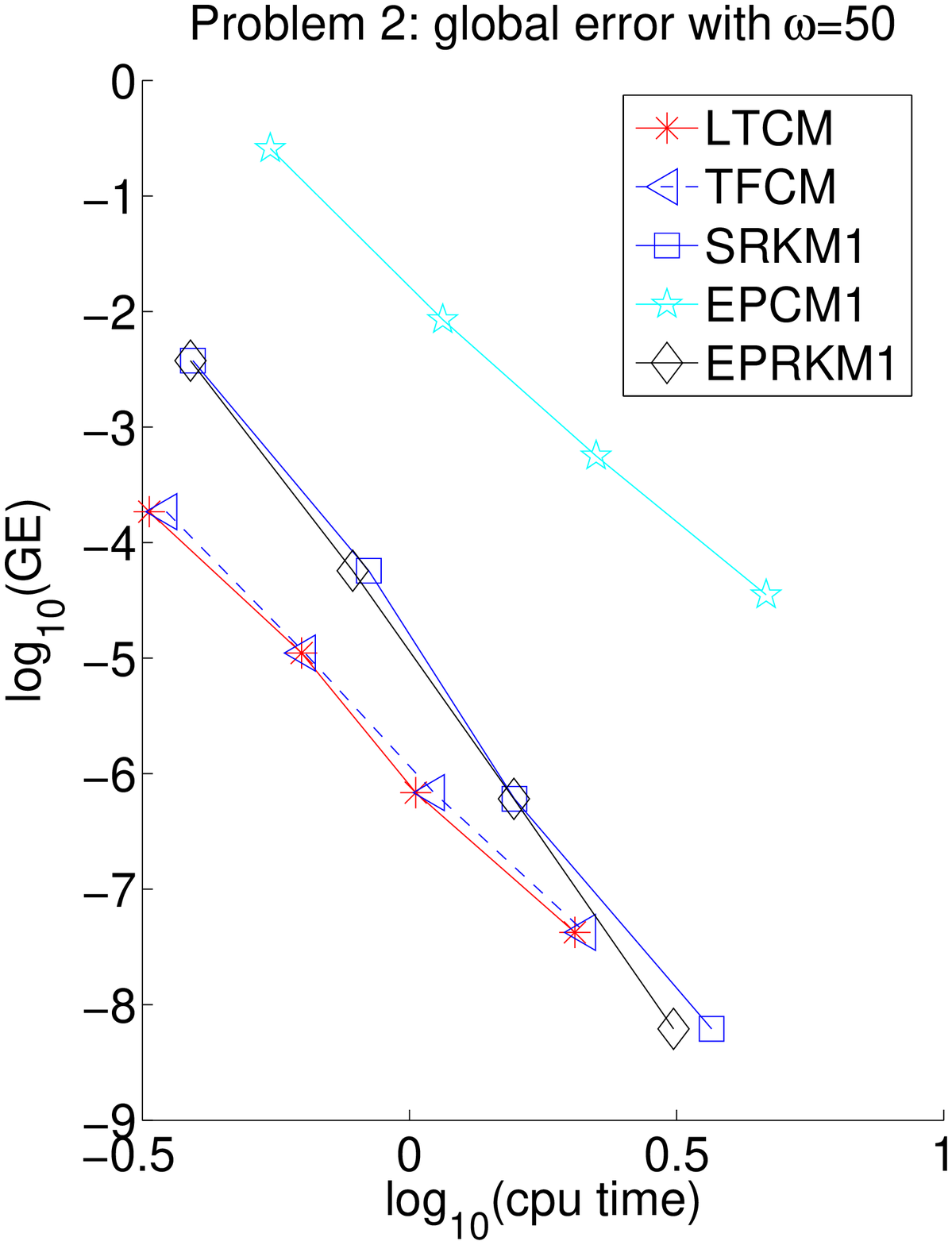} &
\includegraphics[width=5cm,height=6cm]{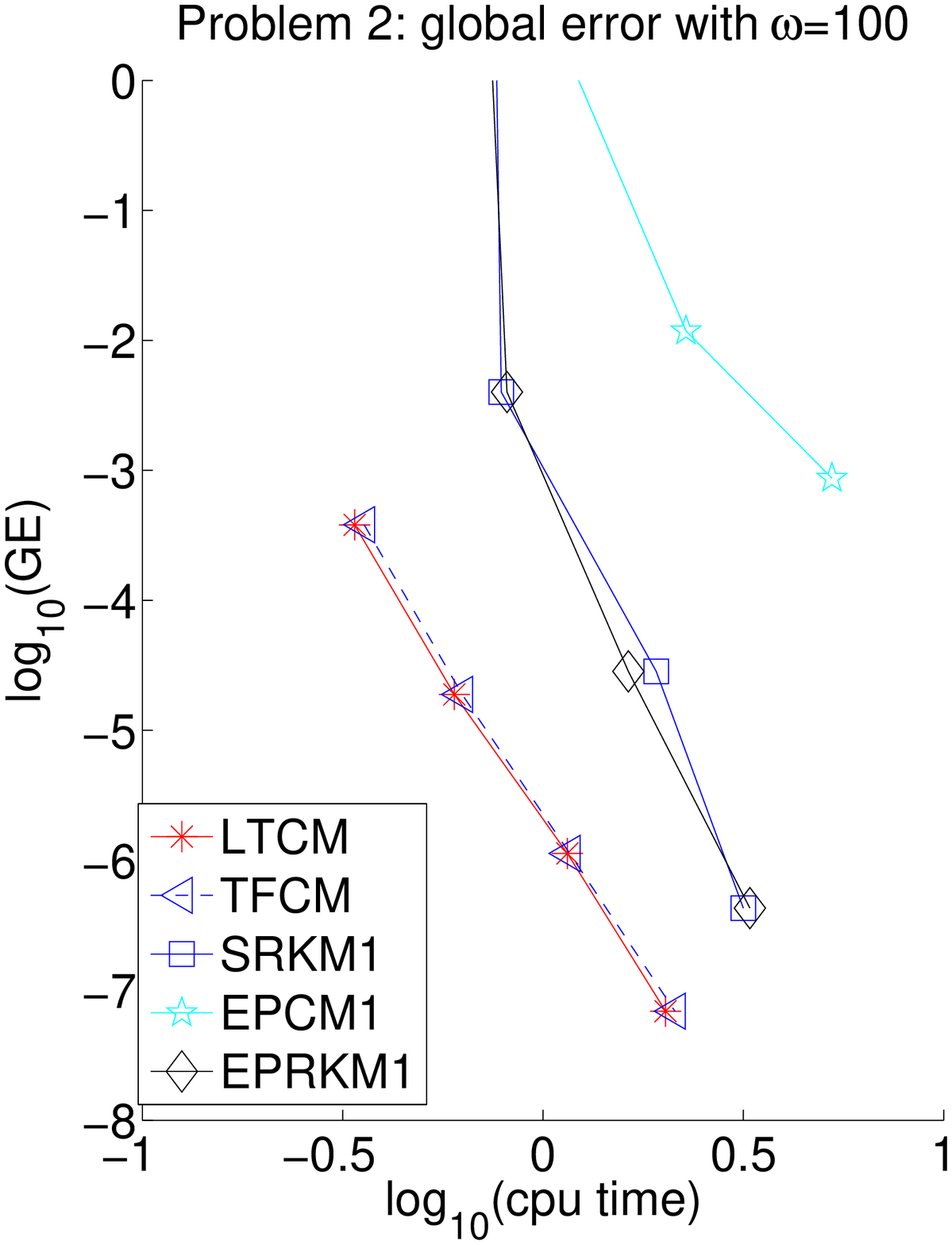}\\
{\small (i)} & {\small (ii)} \\%
\includegraphics[width=5cm,height=6cm]{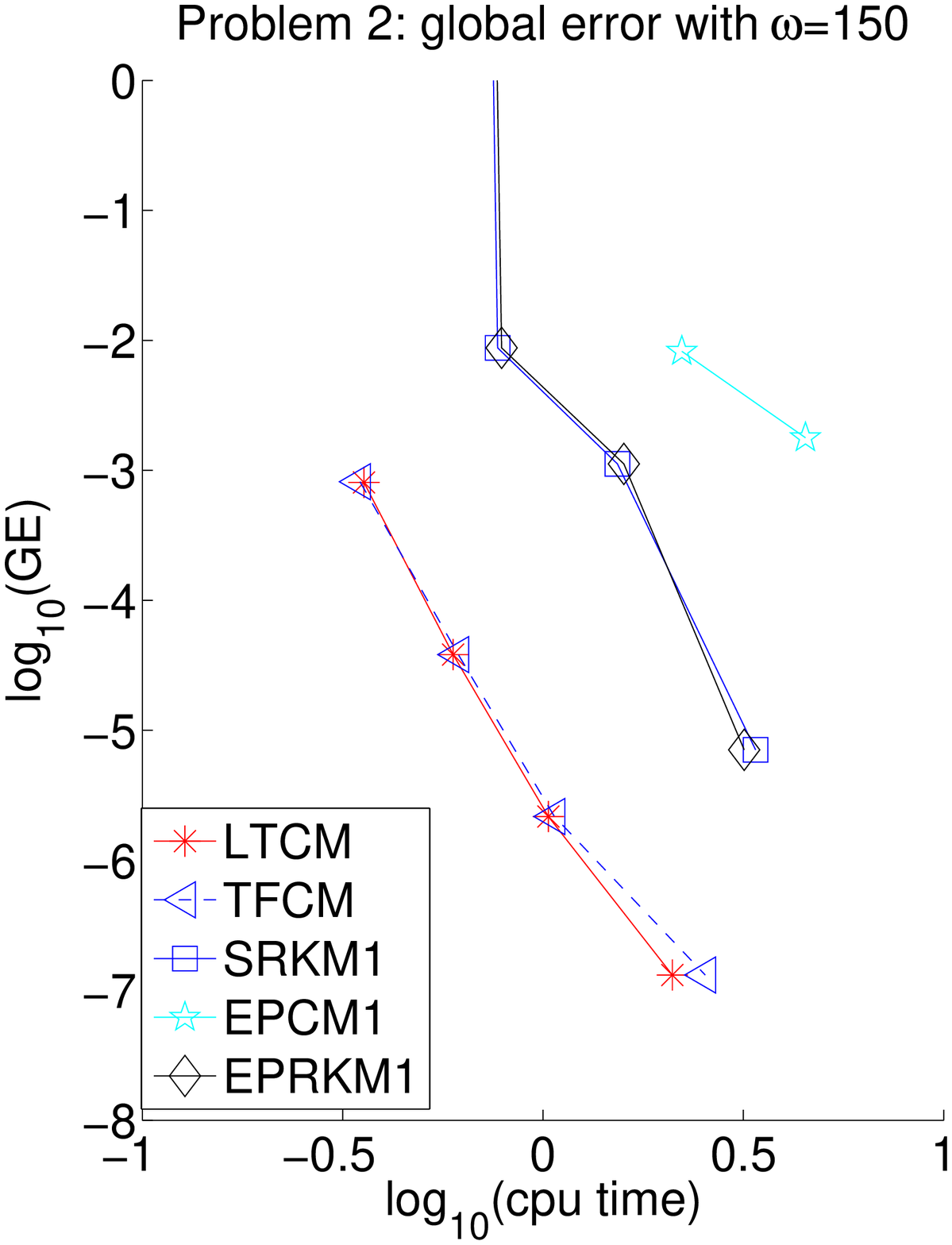} & \includegraphics[width=5cm,height=6cm]{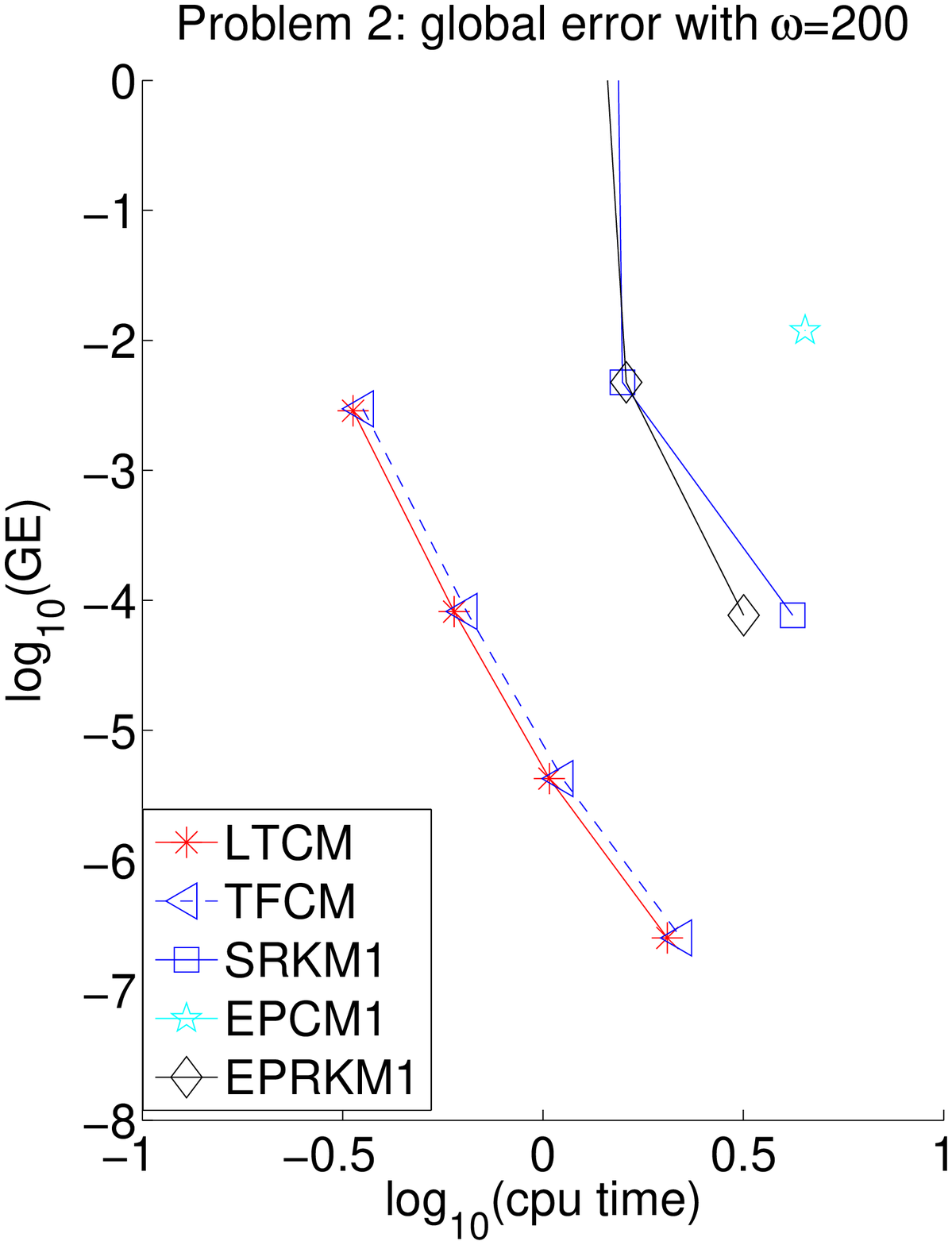}\\
{\small (iii)} & {\small (iv)}%
\end{tabular}
\caption{Results for Problem 2.  The logarithm of the global error
($GE$) over the integration interval against the logarithm of CPU
time. }%
\label{fig:problem2-1}%
\end{figure}

\begin{figure}[ptbh]
\centering\tabcolsep=4mm
\begin{tabular}
[c]{cccc}%
\includegraphics[width=5cm,height=6cm]{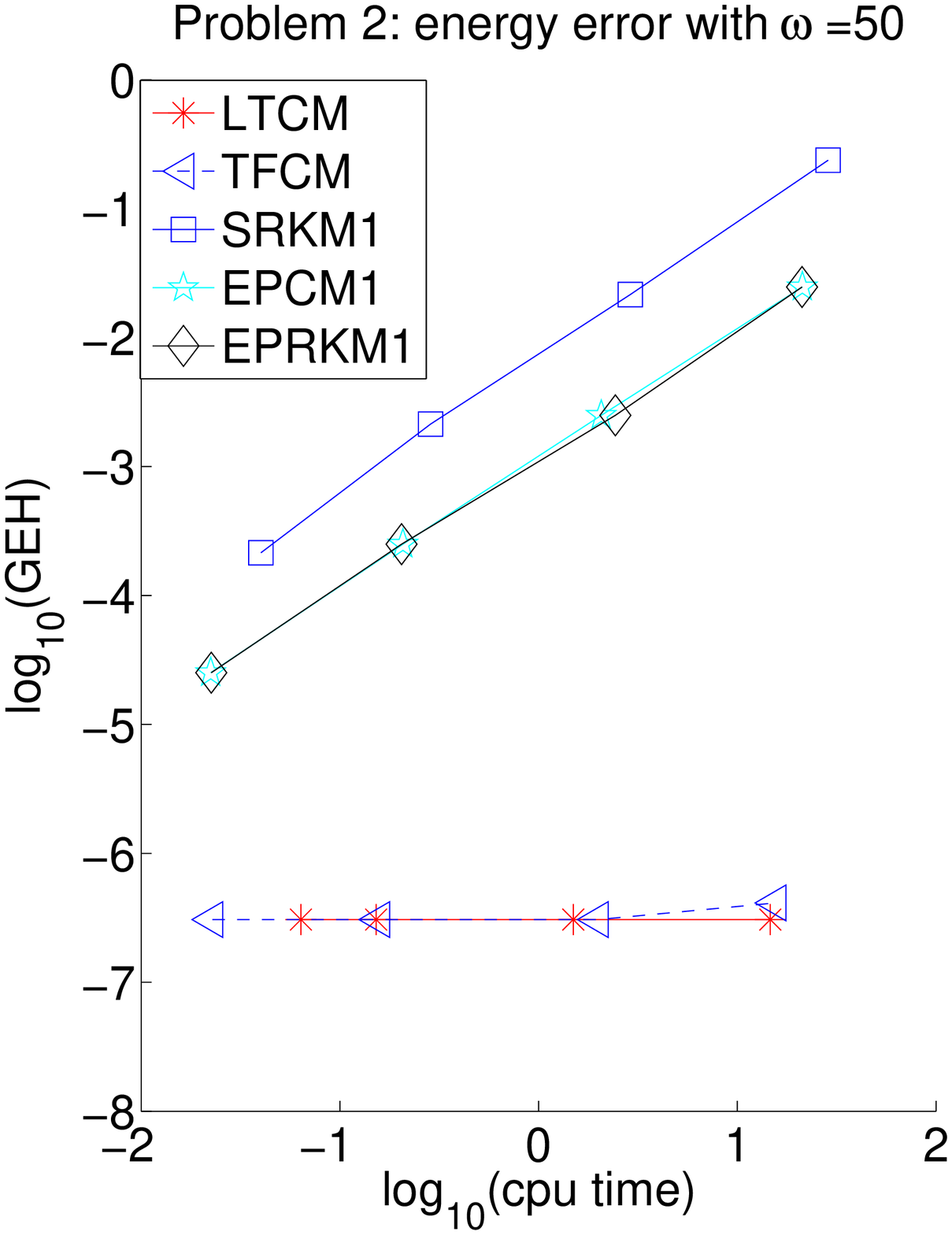} &
\includegraphics[width=5cm,height=6cm]{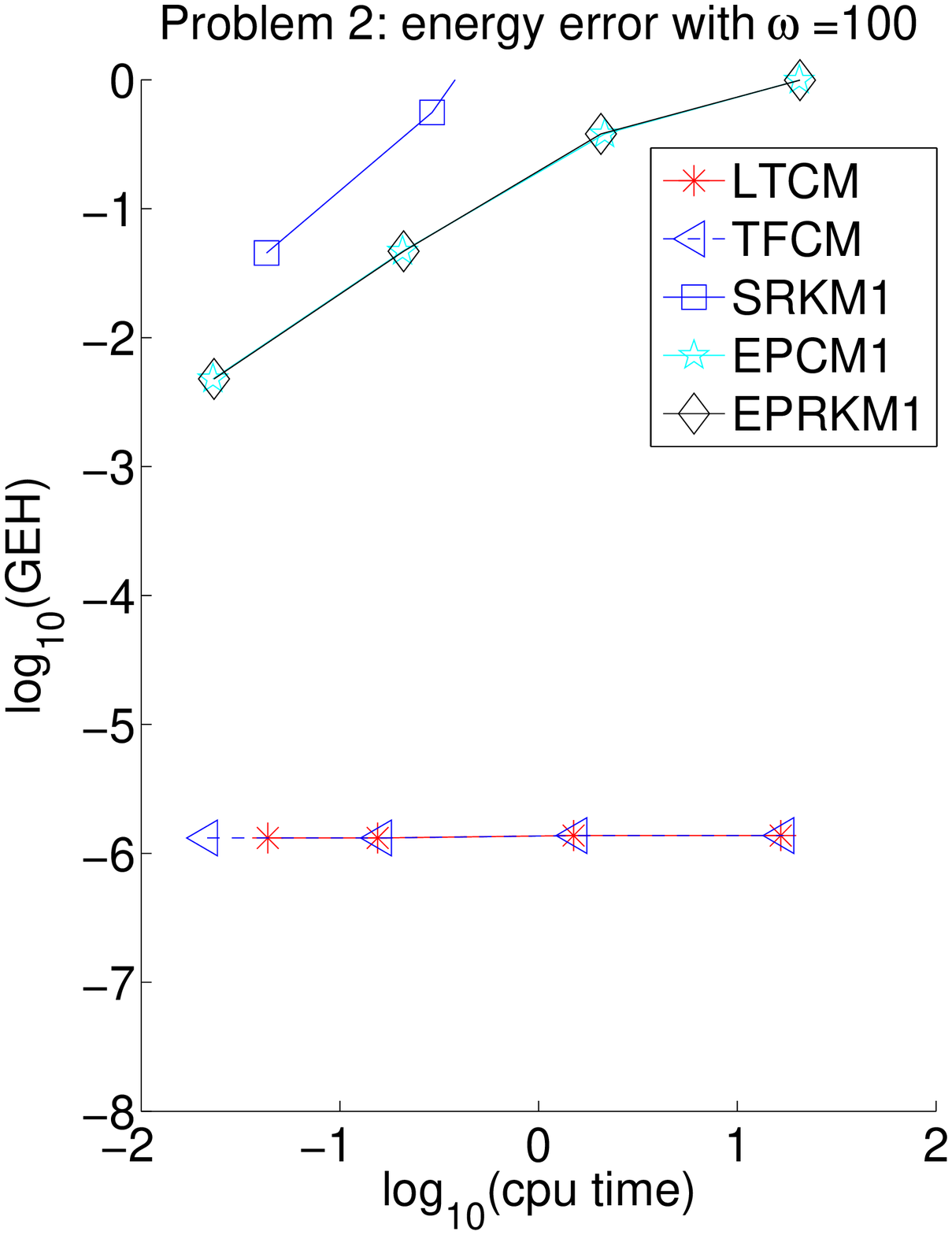}\\
{\small (i)} & {\small (ii)} \\%
\includegraphics[width=5cm,height=6cm]{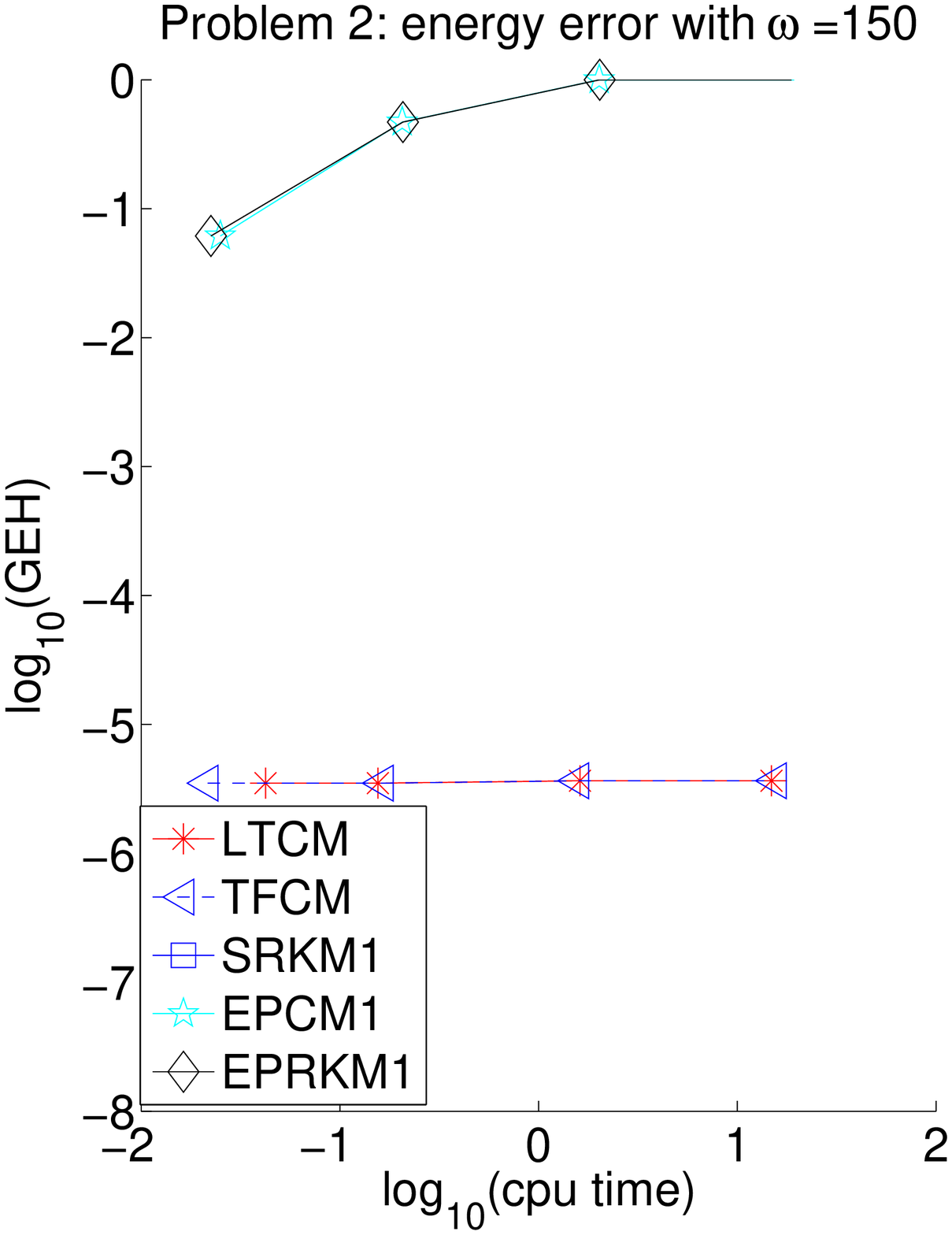} & \includegraphics[width=5cm,height=6cm]{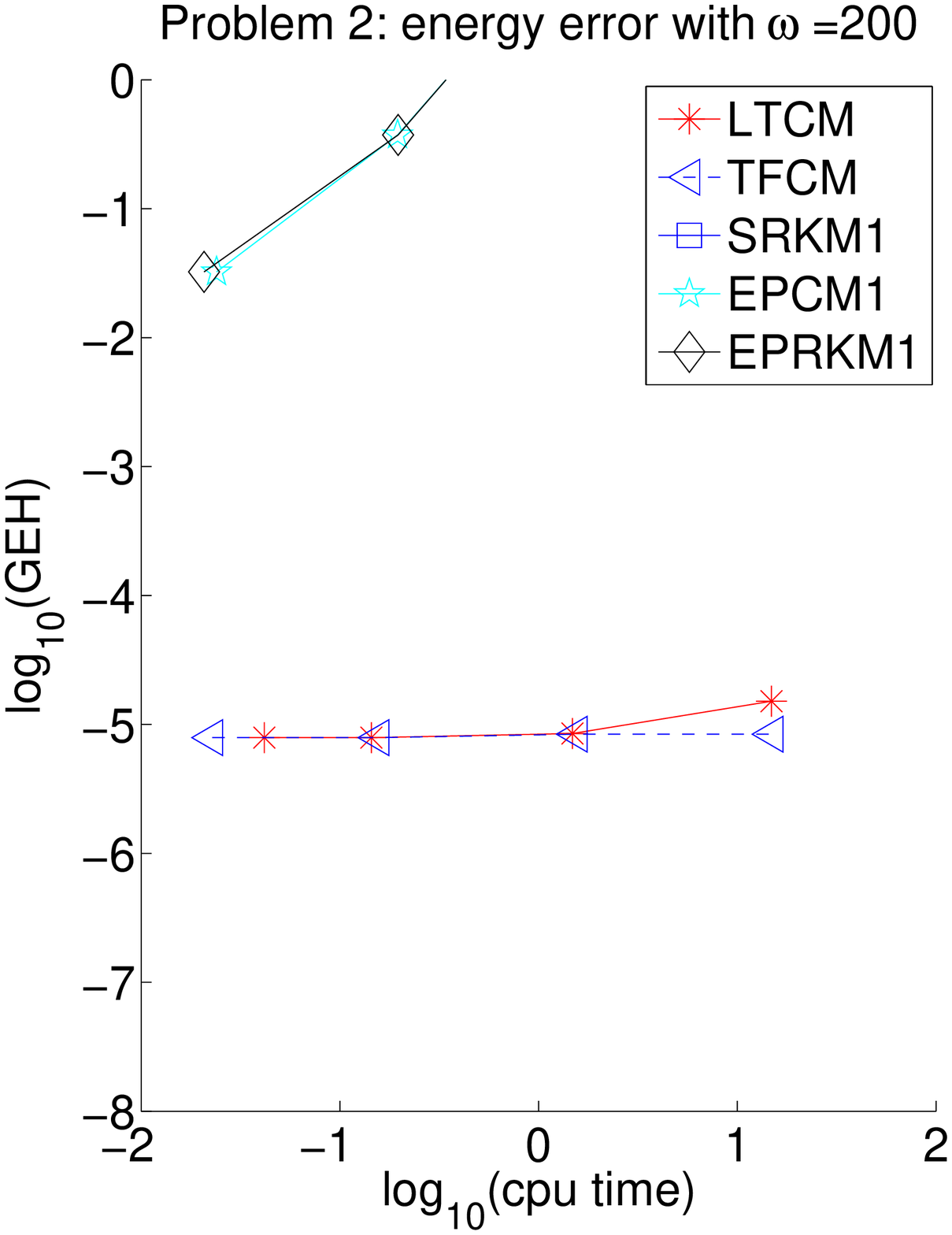}\\
{\small (iii)} & {\small (iv)}%
\end{tabular}
\caption{Results for Problem 2. The logarithm of the maximum global
error of Hamiltonian energy ($GEH$)  against the logarithm of
CPU time.}%
\label{fig:problem2-2}%
\end{figure}

 Then we  integrate the
system   in the interval $[0,50]$ with $\omega=50,100,150,200$ and
the stepsizes $h=1/(20\times{2^{j}}),\ j=1,2,3,4.$ The global errors
are shown in Fig. \ref{fig:problem2-2}.  Finally we integrate this
problem with a fixed stepsize  $h=1/100$  in the interval
$[0,t_{\mathrm{end}}]$ with $t_{\mathrm{end}}=1, 10, 100, 1000.$ The
maximum global errors of Hamiltonian energy  are presented in Fig.
\ref{fig:problem2-2}.
 Here it is noted that some results  are too large,
 thus we do not plot the corresponding points in Figs. \ref{fig:problem2-1}-\ref{fig:problem2-2}.
Similar situation occurs in the next two problems.

\begin{table}$$
\begin{array}{|c|c|c|c|c|c|}
\hline
\text{Methods} &tol=1.0e-006  &tol=1.0e-008   &tol=1.0e-010 &tol=1.0e-012    \\
\hline
\text{LTCM ($\omega=100$)} &  7.1570 &   9.7010  &  9.6435 &  12.2449 \cr
\text{LTCM ($\omega=200$)} &      7.5169 &   10.0160 &    9.2135  &  11.1672\cr
\hline
\text{TFCM ($\omega=100$)} & 7.6434  & 10.3224  & 10.3341  & 12.7998 \cr
\text{TFCM ($\omega=200$)} & 7.8861 &   11.1322 &   10.0578 &   12.3621 \cr
\hline\text{SRKM1 ($\omega=100$)} &32.0491  & 39.4922 &  48.5822  & 57.0720\cr
\text{SRKM1 ($\omega=200$)} &58.2410  &  70.5585  &  86.1757 &   99.6403\cr
\hline\text{EPCM1 ($\omega=100$)}  &   50.8899 &  70.5920  & 87.9782  &102.9839\cr
\text{EPCM1 ($\omega=200$)}  &  121.2714  & 149.7104 &  189.4323  & 220.1096\cr
\hline \text{EPRKM1 ($\omega=100$)}  & 31.0881 &  39.0050 &  47.6389  & 56.4456\cr
\text{EPRKM1 ($\omega=200$)}  & 55.2205  &  68.8459 &   82.5919 &   98.5277\cr
 \hline
\end{array}
$$
\caption{Results for Problem 2: The total CPU time (s) of iterations
for different error tolerances (tol).} \label{pro2-NEW tab-1}
\end{table}

\begin{figure}[ptbh]
\centering\tabcolsep=4mm
\begin{tabular}
[c]{cc}%
\includegraphics[width=5cm,height=6cm]{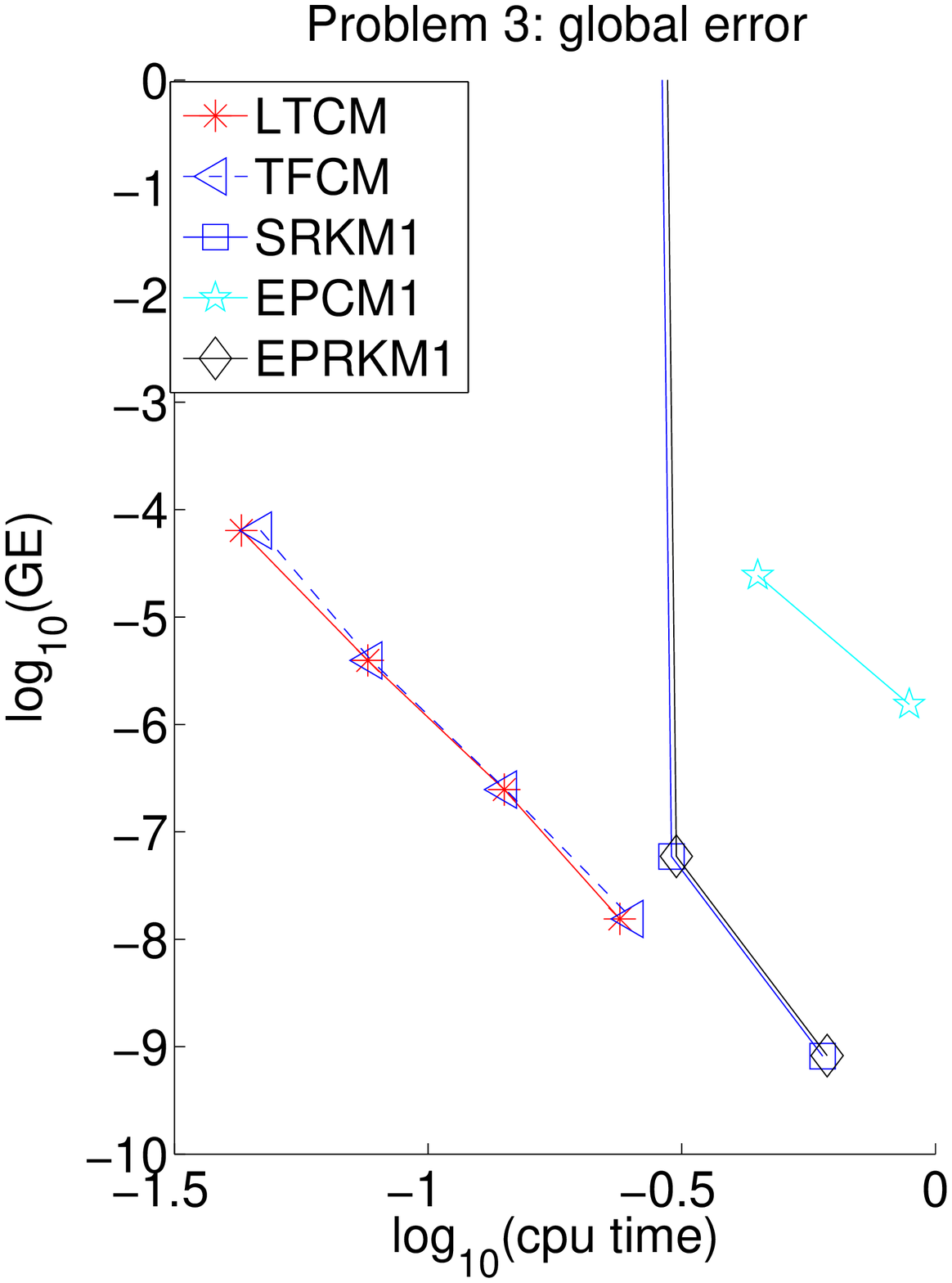} & \includegraphics[width=5cm,height=6cm]{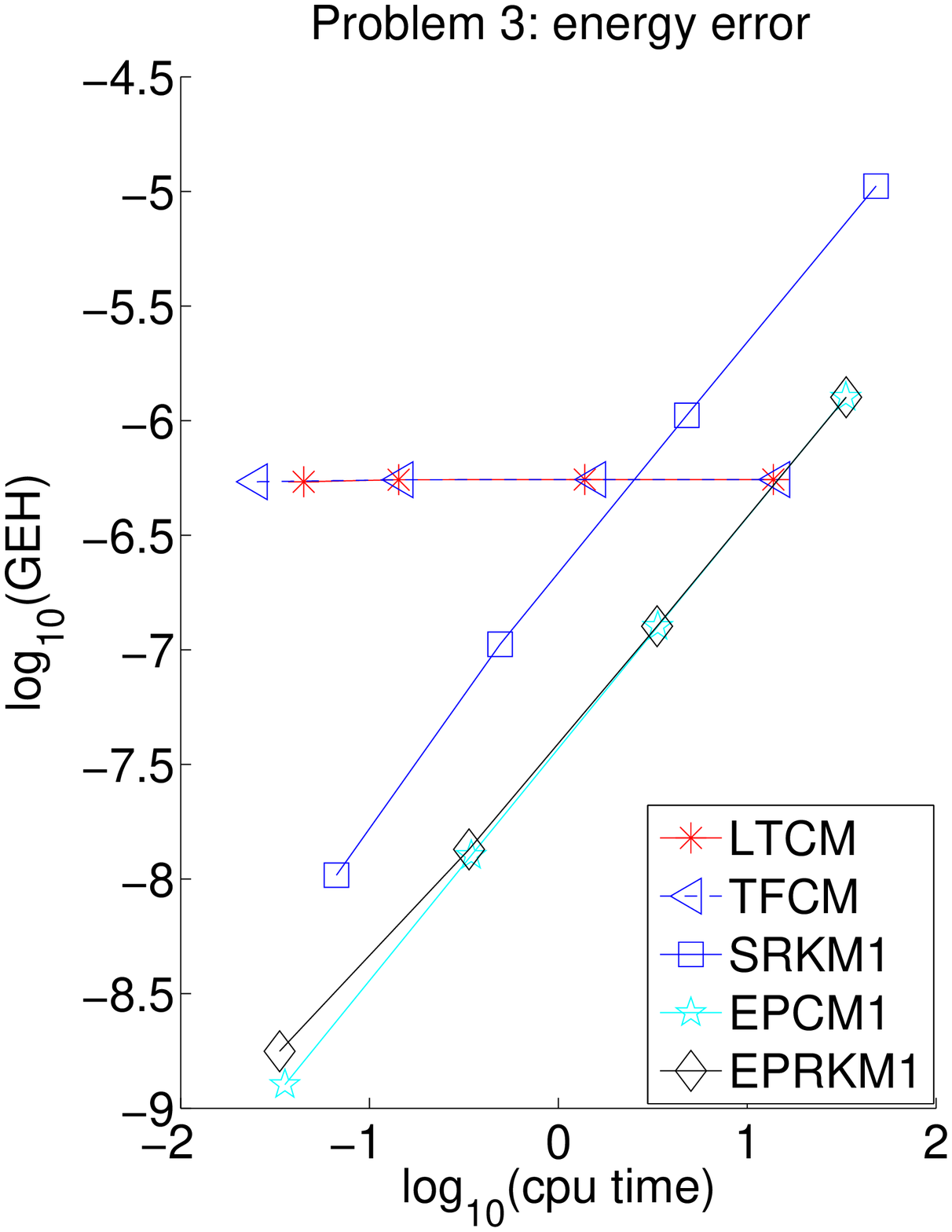}\\
{\small (i)} & {\small (ii)}%
\end{tabular}
\caption{Results for Problem 3. (i): The logarithm of the global
error ($GE$) over the integration interval against the logarithm of
CPU time. (ii):\ The logarithm of the   maximum global error of
Hamiltonian energy ($GEH$)  against the logarithm of
CPU time.}%
\label{fig:problem3}%
\end{figure}

\vskip1mm\noindent\textbf{Problem 3.}
 Consider the
nonlinear Klein-Gordon equation \cite{jimenez1990}
\[
\begin{array}
[c]{l}%
\frac{\partial^{2}u}{\partial t^{2}}-\frac{\partial^{2}u}{\partial
x^{2}}=-u^3-u,\ \ \ 0<x<L,\ \ t>0,\\[0.3cm]%
u(x,0)=A(1+\cos(\frac{2\pi}{L}x)),\ \
 u_{t}(x,0)=0,\ \
u(0,t)=u(L,t),
\end{array}
\]
where $L=1.28$, $A=0.9$. Carrying out a semi-discretization on the
spatial variable by using second-order symmetric differences yields
\begin{eqnarray*}
\begin{array}{l}
\frac{d^2U}{dt^2}+
MU=F(U),\ \ \ 0<t\leq t_{\mathrm{end}},\\
\end{array}
\end{eqnarray*}
where $U(t)=\big(u_1(t),\cdots,u_N(t)\big)^{\intercal}$ with
$u_i(t)\approx u(x_i,t),\ i=1,\cdots,N$,
\begin{eqnarray*}
M=\frac{1}{\Delta x^2}\left(
\begin{array}
[c]{ccccc}%
2 &-1 && &-1\\
-1 &2  & -1&  &  \\
 &\ddots&\ddots&\ddots& \\
&&-1 &2  & -1\\
 -1 & & &-1&2  \\
\end{array}
\right)_{N\times N}
\end{eqnarray*}
with $\Delta x= L/N$,  $x_i = i\Delta x,$ $
F(U)=\big(-u_1^3-u_1,\cdots,-u_N^3-u_N\big)^{\intercal}$ and $N=32$.
The corresponding
  Hamiltonian of this system is
\begin{equation*}
H(U',U)=\frac{1}{2}U'^{\intercal}U'+\frac{1}{2}U^{\intercal}MU+\frac{1}{2}u^2_{1}+\frac{1}{4}u^4_{1}+\ldots+
\frac{1}{2}u^2_{N}+\frac{1}{4}u^4_{N}.
\end{equation*}
Here we choose $N=32$. The problem is solved in the interval $[0,
500]$ with the stepsize $h=\frac{1}{100}$ to show the convergence
rate of iterations. See Table \ref{pro3-NEW tab} for the total
 CPU time of iterations for
different error tolerances.  Then we   solve this problem  in
$[0,20]$ with stepsizes
  $h=1/(3\times2^{j}),\ j=1,2,3,4.$ Fig.
\ref{fig:problem3} (i) shows the global errors. Finally this problem
is integrated  with a fixed stepsize  $h=0.002$ in the interval
$[0,t_{\mathrm{end}}]$ with $t_{\mathrm{end}}=1, 10, 100, 1000.$ The
maximum global errors of Hamiltonian energy   are presented in Fig.
\ref{fig:problem3} (ii).

\begin{table}$$
\begin{array}{|c|c|c|c|c|c|}
\hline
\text{Methods} &tol=1.0e-006  &tol=1.0e-008   &tol=1.0e-010 &tol=1.0e-012    \\
\hline
\text{LTCM} &  5.9325&    7.9263  &  8.1816 &  10.0602 \cr
\text{TFCM} & 6.5318 &   8.7008  & 8.8934 &  10.7489 \cr
\text{SRKM1} &   24.1600&   29.4173 &  34.5310 &  39.5161\cr
\text{EPCM1}  &37.2757 &  46.4011 &  53.1403 &  66.2339\cr
\text{EPRKM1}  & 22.6571 &  27.8341 &  33.5435  & 39.4533\cr
 \hline
\end{array}
$$
\caption{Results for Problem 3: The total CPU time (s) of iterations
for different error tolerances (tol).} \label{pro3-NEW tab}
\end{table}

\noindent\vskip3mm \noindent\textbf{Problem 4.} Consider the wave
equation
\begin{eqnarray*}\begin{array}{ll} \frac{\partial^2u}{\partial
t^2}-a(x)\frac{\partial^2u}{\partial
x^2}+92u=f(t,x,u),\ \ \ 0<x<1,\ \ t>0,\\\\
u(0,t)=0,\ \ \ u(1,t)=0,\ \ \ u(x,0)=a(x),\ \ \ u_t(x,0)=0,
\end{array}
\end{eqnarray*}
with $ a(x) = 4x(1-x),\
f(t,x,u)=u^5-a^2(x)u^3+\frac{a^5(x)}{4}\sin^2(20t)\cos(10t).$ The
exact solution is  $u(x,t) = a(x) \cos(10t).$

\begin{figure}[ptbh]
\centering\tabcolsep=4mm
\begin{tabular}
[c]{c}%
\includegraphics[width=5cm,height=6cm]{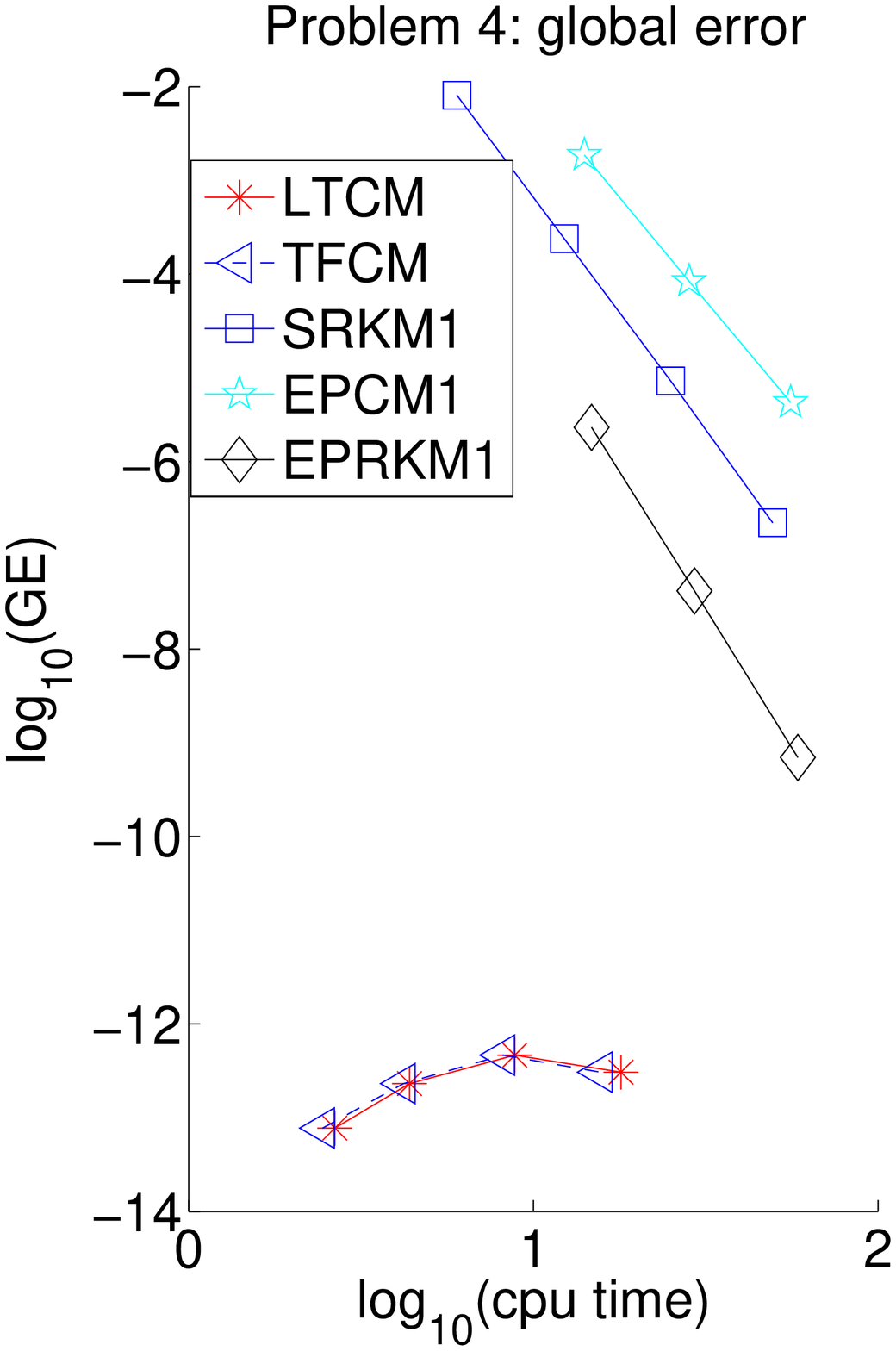}
\end{tabular}
\caption{Results for Problem 4: The logarithm of the global error
($GE$) over the integration interval against the logarithm of
CPU time.}%
\label{fig:problem4}%
\end{figure}

Using semi-discretization on the spatial variable with second-order
symmetric differences, we obtain
\begin{eqnarray*}
\begin{array}{ll}
\frac{d^2U}{dt^2}+MU=F(t,U),\
U(0)=\big(a(x_1),\cdots,a(x_{N-1})\big)^{\intercal},\ U'(0)={\bf0},
\ 0<t\leq t_{\mathrm{end}},
\end{array}\label{pro4}
\end{eqnarray*}
where $U(t)=\big(u_{1}(t),\ldots,u_{N-1}(t)\big)^{\intercal}$ with
$u_{i}(t)\approx u(x_{i},t)$, $x_i = i\Delta x$, $\Delta x= 1/N$,
$i=1,\ldots,N-1,$
\begin{eqnarray*}
M=92I_{N-1}+\frac{1}{\Delta x^2}\left(
\begin{array}
[c]{ccccc}%
2a(x_1)&-a(x_1) && &\\
-a(x_2) &2a(x_2) & -a(x_2)&  &  \\
 &\ddots&\ddots&\ddots& \\
&&-a(x_{N-2}) &2a(x_{N-2}) & -a(x_{N-2})\\
  & & &-a(x_{N-1})&2a(x_{N-1}) \\
\end{array}
\right),
\end{eqnarray*}
 and
\begin{eqnarray*}
F(t,U)=\big( f(t,x_1,u_1), \cdots,
f(t,x_{N-1},u_{N-1})\big)^{\intercal}.
\end{eqnarray*}
The problem is solved in the interval $[0, 100]$ with the stepsize
$h=\frac{1}{40}$ to show the convergence rate of iterations. See
Table \ref{pro4-NEW tab} for the total
 CPU time of iterations  for different error tolerances.  Then
system is integrated in the interval $\lbrack0,100]$ with $N=40$ and
$h=1/2^j,\ j=5,6,7,8.$ The   global errors are shown in Fig.
\ref{fig:problem4}.

 It follows from the numerical results that our method LTCM is very promising
as compared with the classical methods  SRKM1, EPCM1 and EPRKM1.
Although LTCM has a similar performance as TFCM in preserving the
solution and the energy, it has a better convergence rate of
iterations.
\begin{table}$$
\begin{array}{|c|c|c|c|c|c|}
\hline
\text{Methods} &tol=1.0e-006  &tol=1.0e-008   &tol=1.0e-010 &tol=1.0e-012    \\
\hline
\text{LTCM} &  1.8980   &   1.8737 &     2.1212  &    2.3196 \cr
\text{TFCM} &  1.9213   &   1.9345   &   2.2227   &   2.3736 \cr
\text{SRKM1} & 13.8634  &   16.6963  &   19.0854   &  22.6142\cr
\text{EPCM1}  &23.5110   &  28.1288   &  32.2263   &  36.8443\cr
\text{EPRKM1}  &13.5526  &   17.2289  &   18.8744   &  23.0066\cr
 \hline
\end{array}
$$
\caption{Results for Problem 4: The total CPU time (s) of iterations
for different error tolerances (tol).} \label{pro4-NEW tab}
\end{table}
\section{Conclusions and discussions}
\label{sec:conclusion} In this paper we have investigated a kind of
trigonometric collocation methods based on Lagrange basis
polynomials, the variation-of-constants formula and the idea of
collocation methods for solving
 multi-frequency oscillatory second-order differential
equations   \eqref{prob} efficiently. It has been shown that the
convergent condition of these   trigonometric collocation methods is
independent of $\norm{M}$, which is very important and crucial for
solving highly oscillatory systems. The numerical experiments with
some model problems show that the our   method derived in this paper
has remarkable efficiency in comparison with some existing methods
in the
 literature.
\section*{Acknowledgements}
The authors are sincerely thankful to two anonymous reviewers for
their valuable suggestions, which help improve the presentation of
the manuscript significantly.

\end{document}